\documentclass[12pt,a4paper]{article}
\usepackage[english]{babel}
\usepackage{bm}
\usepackage{amsmath,amsthm,amssymb}
\usepackage{amsthm}
\usepackage{cases}
\usepackage{latexsym}
\usepackage{comment}
\usepackage{enumerate}
\usepackage{xcolor}
\usepackage{courier}        
\usepackage{graphicx}
\usepackage{dsfont}

\theoremstyle{plain}
\newtheorem{thm}{Theorem}
\newtheorem{lem}[thm]{Lemma}

\newtheorem{rem}{Remark}
\newtheorem{cor}{Corollary}

\theoremstyle{definition}
\newtheorem{defn}{Definition}[section]

\theoremstyle{remark}

\begin{document}
\title{Quantitative propagation of chaos for 2D stochastic vortex model on the whole space under moderate interactions}

\author{Alexandre B. de Souza \footnote{Departamento de Matem\'atica, Universidade Estadual de Campinas, Brazil. \texttt{a265040@dac.unicamp.br}.} \and 
}

\date{}
\maketitle

\begin{abstract}
We derive the stochastic 2D vortex model on the whole Euclidean space from moderately interacting particle systems driven by individual and environmental noises, obtaining quantitative estimates in the sense of the entropy and energy functionals. The main novelties lie in combining the control of the Fisher information of the particle system with the Ladyzhenskaya and Donsker–Varadhan inequalities, as well as localization techniques within the probabilistic data setting, to address the nonlinearity and quadratic variation arising from Itô's formula.
Moreover, we construct a suitable solution for the limiting process. 
\end{abstract}

\vspace{0.3cm} \noindent {\bf MSC2010 subject classification:} 35R60,  35Q30, 82C22

\section{Introduction}
In this work, we  establish a quantitative derivation of the stochastic 2D
vortex model in the whole Euclidean space, which is a particular case of the following stochastic Fokker–Planck equation
\begin{align}\label{SPDE_Ito}
\mathrm{d}\rho_t =  \Delta \rho_t \,dt + \frac{1}{2}\nabla^2 \rho_t (\sigma \sigma^\top)_t \, dt - \nabla \cdot( \rho_t (K\ast \rho_t)) \, dt  - \nabla \rho_t \cdot \sigma_t\, d  B_t,
\end{align}
from a stochastic moderately interacting particle system given by     
\begin{align}\label{particles}
d  X_t^{i,N} &=  \frac{1}{N}\sum_{k=1}^{N} \left(K \ast V^{N}\right)\left(X_{t}^{i,N} -X_{t}^{k,N}\right) \, dt + \sqrt{2} \,d W_t^{i,N} + \sigma_t\, dB_t,  \quad 
\end{align}
	\noindent where $W_{t}^{i,N}$ and $B_{t}$ are  independent standard $\mathbb{R}^d$-valued Brownian motions, defined on a filtered probability space $(\Omega,\mathcal{F},(\mathcal{F}_t)_{t\geq0},\mathbb{P})$, $K$ is a singular kernel,  $V^{N}$ is a suitable scaling and $\sigma$ is a diffusion matrix.

To achieve it, we study the asymptotic behavior, as $N \to \infty$ of the \emph{empirical measure} defined by
\begin{align*}
    S_{t}^{N}\doteq\frac{1}{N} \sum_{i=1}^{N}\delta_{X_{t}^{i,N}}, \qquad t\geqslant 0
\end{align*}
by introducing the \emph{mollified empirical measure} 
\begin{align*}
    \rho_t^{N}\doteq V^{N} \ast S_t^{N} = \int_{\mathbb{R}^d} V^N(\cdot-y) \,S_t^N(dy),
\end{align*}
which allows us to apply the entropy and energy functionals directly to the densities $\rho^N$ and $\rho$. 

Particle systems driven by environmental noises appear in several contexts, such as in statistical mechanics, population dynamics, as well as in numerical analysis, see, e.g., \cite{Carmona}, \cite{Carmona2}, \cite{Jabin2}, \cite{Live}, and \cite{Lac}, and quantitative bounds characterize a reduction in the complexity of the system (\ref{particles}), for large values of $N$, see \cite{Jabin_2018}.

To the best of our knowledge, this is the first time that, propagation of chaos, supplemented with pathwise quantitative bounds, in the sense of the relative entropy functional, have been derived for the 2D stochastic vortex model in the whole Euclidean space via moderately interacting particles. This extends previously known results for bounded kernels on $\mathbb{R}^d$ \cite{Niko} and for sub-Coulombian kernels under periodic boundary conditions \cite{Ale}.
In addition, we obtain energy estimates for the difference between the regularized empirical measure and the solution of (\ref{SPDE_Ito}), extending previous results by \cite{Josue} and \cite{Pisa}.

Our main novelty lies in the key technical tools we combine. Using Itô formula, we compute the evolution equations for the entropy and energy functionals associated with the regularized empirical measure and the solution of the aforementioned model. By combining the control of the Fisher information of the particle system with the Ladyzhenskaya and Donsker–Varadhan inequalities, along with localization techniques within the probabilistic data setting, to deal with the nonlinearity and quadratic variation, we are able to close the argument with a nonlinear Grönwall inequality. Moreover, we prove the existence of a suitable solution to equation~\eqref{SPDE_Ito}.

\subsection*{Related works}
Quantitative estimates  for singular interactions have been substantially developed in recent years.
The relative entropy method for establishing quantitative propagation of chaos with singular interaction kernels in $W^{-1,\infty}$ was first introduced in \cite{Jabin_2018}, providing a framework that includes the 2D viscous vortex model. 
Alternatively, we refer to \cite{Serfaty}, which establishes quantitative propagation of chaos for McKean–Vlasov equations with singular interaction kernels via the modulated energy method.
 In contrast to the relative entropy approach, which operates at the level of the joint law of the particle system, the modulated energy method is formulated directly in terms of the empirical measure associated with the particles.

Since then, the aforementioned works have paved the way for establishing quantitative propagation of chaos results, under singular interactions. In \cite{Bres} and \cite{Bres1}, the strategies given in \cite{Jabin_2018} and \cite{Serfaty}, were adapted to obtain a quantitative derivation for the  Patlak–Keller–Segel model, in optimal subcritical regimes. 
Additionally, \cite{Gui}, by means of a logarithmic Sobolev inequality, derived uniform in time propagation of chaos for the 2D vortex model,   adapting the strategy in \cite{Jabin_2018}. Alternatively, with sharp estimates, global in time derivation of that model, based on 
relaxation estimates for the limiting equation and inequalities of Riesz modulated energy, was derived by \cite{Cho2}, see also  \cite{Cho}, \cite{Nguyen} and \cite{Rosey1}, for the modulated energy method.
Regarding sub-Coulombian kernels,
\cite{Carri} obtained the derivation of the mean-field approximation for Landau-like equations, while \cite{Rosey} obtained global-in-time mean-field convergence for gradient (conservative) diffusive flows.
Moreover, quantitative entropy estimates  for 
systems with individual and common noise, have been obtained, as shown in \cite{Shao} for the 2D vortex model, on torus, \cite{Chen2} for the Hegselmann-Krause model.

More recently, quantitative propagation of chaos for the 2D vortex model on the whole space was obtained in \cite{Feng, Feng1}, while the corresponding result for the 2D log-gas was established in \cite{Feng2}. In these works, the authors exploit Li–Yau-type estimates and Hamilton-type heat kernel estimates to derive suitable bounds for the limiting solution and its derivatives up to second order, which allows them to extend the strategy given in \cite{Jabin_2018}, to $\mathbb{R}^d$.

Quantitative estimates in the whole space were derived by \cite{Pisa} and \cite{Pisa1} in the context of moderately interacting particles\footnote{The moderately interacting particle model was introduced by Oelschläger in \cite{Oelschlager84, Oelschlager85, Oelschlager87}.} without common noises and with singular kernels (including the Biot-Savart, Keller-Segel, and Dirac delta measure kernels). Their derivation is based on a semigroup approach introduced in \cite{FlandoliLeimbachOlivera} and has motivated further advancements, for example, in \cite{ric, correa, Hao, Josue1, Simon2}.
Under periodic boundary conditions, \cite{Josue} used the Krilov $L^p$ theory for SPDE's, to extend the results given in \cite{Pisa}, in the presence of environmental noises, as the semigroup approach is not available in this setting.
It is important to note that these works do not rely on the relative entropy method, but instead establish quantitative estimates in strong functional topologies.

The connection  between the time evolution of the relative entropy, of the joint
distribution, and moderate interacting regime, on the whole space was carried out by  \cite{Chen}, based on a combination between the relative entropy and the regularised $L^{2}$-estimate in \cite{Oelschlager87}, deriving a propagation of chaos result for the viscous porous medium equation. However, it was in \cite{Ale} that equation \eqref{SPDE_Ito} was derived on the torus for singular kernels and common noise in the moderately interacting regime, by applying the relative entropy directly to the mollified empirical measure. 

In this work, we extend the quantitative entropy estimates given in \cite{Ale} under periodic boundary conditions and in \cite{Niko} for bounded kernels to the whole space for the stochastic 2D vortex model. Additionaly, we obtain energy estimates for that model, extending the results in \cite{Josue} and \cite{Pisa}.
In contrast to the aforementioned works, which establish quantitative entropy estimates at the level of the joint law, our application of the relative entropy method provides pathwise estimates, yielding bounds directly at the trajectory level of the particle system.

\subsection*{Outline of the article}
In order to present our results, this work is organized as follows. In the next three subsections, \ref{notations}, \ref{assumpions}, and \ref{solution}, we introduce the notations and assumptions used throughout the text, setting the stage for the subsequent
discussions, as well as the notion of solution to (\ref{SPDE_Ito}). In the subsection \ref{mainr}, we present our main results. The sections \ref{Proof} and \ref{Proof2} are dedicated to the proofs of our main results. Finally, the Appendix provides the existence of a suitable solution for the limiting equation, as well as some inequalities used throughout the proofs.

\subsection{Notations} \label{notations}
	For \( d \geq 1 \), 
    the space of probability density functions on \( \mathbb{R}^d \) is denoted by \( \mathcal{P}(\mathbb{R}^d) \).	For a normed vector space \(U\), we denote its norm by \(\|\cdot\|_U\), except for \(L^p\) spaces, where we write \(\|\cdot\|_p\), for $p\geq 1$.
    
    For $\gamma \in (0,1]$ the H\"older space on $\mathbb{R}^d$ is given by $$C^{\gamma}=C^{\gamma}(\mathbb{R}^d)\doteq \left\{f:\mathbb{R}^d \to \mathbb{R}^{e} \mid\|f\|_{\gamma}\doteq \|f\|_{\infty} + \sup_{x,y\in \mathbb{R}^d}\frac{|f(x)-f(y)|}{|x-y|^{\gamma}}<\infty\right\}.$$
The Kantorovich-Rubinstein metric reads, for any two probability measures $\mu$ and $\nu$ on $\mathbb{R}^d$,			\begin{equation}\label{eq:defWasserstein}				\|\mu - \nu \|_{0} \doteq \sup \left\{ \int_{\mathbb{R}^d} \phi \, d(\mu-\nu) \, ; \|\phi\|_\infty,  \|\phi\|_{\text{Lip}} \leq 1 \right\} .
		\end{equation}
Let $f$ and $g$ be positive probability density functions on \(\mathbb{R}^d\). The relative entropy (or Kullback–Leibler divergence) of $f$ with respect to $g$ is defined as 
\begin{align*}
    \mathcal{H} \left(f|g\right)\doteq \int_{\mathbb{R}^d}f(x)\ln{\left(\frac{f(x)}{g(x)}\right)} \, dx.
\end{align*}
Also the Fisher information of $f$ with respect to $g$ is given by
\begin{align*}
    \mathcal{I}\left(f|g\right) \doteq \int_{\mathbb{R}^d}f(x)\left|\nabla \ln{\frac{f(x)}{g(x)}}\right|^2 \, dx = \int_{\mathbb{R}^d}\frac{(g(x))^2}{f(x)}\left|\nabla \left(\frac{f(x)}{g(x)}\right)\right|^2 \, dx.
\end{align*}

\subsection{Assumptions}\label{assumpions}
			
\begin{enumerate}
\item[$(\mathbf{A}^V)$] Let  $V^N:\mathbb{R}^d\to \mathbb{R}$ be the scaling given by $V^N(y) \doteq N^{\beta}V(N^{\frac{\beta}{d}}y)$, $\beta \in (0,1)$, with $V \in \mathcal{P}\left(\mathbb{R}^d\right)$ satisfying
\begin{align}
    \left|\nabla V\right| \leq C_{d} V, \label{mollifier inequality}
\end{align}
for $\alpha >d/2$
\begin{align}
    V(y) \leq C_{d,\alpha} (1 + |y|^2)^{-\alpha}, \label{decay mollifier}
\end{align}
and for $\gamma, \beta \in (0,1)$ and $y\neq 0$
\begin{align}
    V(y) \leq C_{d,\gamma,\beta} |y|^{-3d/\beta - \gamma}. \label{decay mollifier12}
\end{align}
\item [$(\mathbf{A}^{ K})$]
It holds for $\gamma \in (0,1/2)$ that 
\begin{align}
    \|K \ast f\|_\gamma \leq C_{\gamma}\left(\|f\|_1 + \|f\|_4\right) \label{Holderbiot}
\end{align}
for all $f \in W^{1,2}$.
\item[$(\mathbf{A}^{\nabla\cdot K})$] 
It holds $\nabla\cdot K =0$ and there exists $K_0 \in L^\infty$, with $\|K_0\|_{\infty}\leq 1/4$, such that 
\begin{align}
    K= \nabla\cdot K_0. \, 
\end{align}

\item[$(\mathbf{A}^{\sigma})$]
	The  coefficient  $\sigma: [0,T] \to \mathbb{R}^{d \times d}$ is measurable and bounded.
 \item[$(\mathbf{A}^{\rho_0})$]
	The initial condition $\rho_0$ is taken such that
    $\rho_0 \in W^{2,1} \cap W^{2,\infty}\cap \mathcal{P}\left(\mathbb{R}^d\right)$, and
    \begin{align}
        \left|\nabla \ln{\rho_0}(x)\right|^2 \leq \widetilde{C}_1(1+|x|^2),
    \end{align}
     \begin{align}
        \left|\nabla^2 \ln{\rho_0}(x)\right| \leq \widetilde{C}_2(1+|x|^2)
    \end{align}
    and
     \begin{align}
       \rho_0(x) \leq \widetilde{C}_3\exp{\left(-\widetilde{C}_3^{-1}|x|^2\right)}. \label{decayi}
    \end{align}
    \item[$(\mathbf{A}^{X_0})$]
The initial particle system $(X_0^{i,N})_{i=1}^N$, $N\in \mathbb{N}$ is a sequence of i.i.d random variables with law $\rho_0$.
	\end{enumerate}  
    
\begin{rem}
    For an example of mollifier  $V \in \mathcal{P}(\mathbb{R}^d)$ satisfying $\left(\mathbf{A}^{V}\right)$, take
\begin{align*}
    V(y)  \doteq \overline{C}_{d}\exp{\left(-\sqrt{1+|y|^2}\right)}
     \end{align*} 
with $\overline{C}_{d} \doteq \left(\int_{\mathbb{R}^d}\exp{\left(-\sqrt{1+|y|^2}\right)}\,dy \right)^{-1}$.
\end{rem}

\begin{rem}
    The Biot-Savart kernel in dimension two is given by
    \begin{align}
        K(x)=\frac{1}{2\pi}\frac{x^\perp}{|x|^2}
    \end{align}
    and then, it is the divergence in the sense of distributions, of an $L^\infty$ matrix $K_0$, namely $K=\nabla \cdot K_0$, with
    \begin{align}
        K_0(x_1,x_2) = -\left(\frac{1}{2\pi}\arctan{\left(\frac{x_1}{x_2}\right)}\right)I.
    \end{align}
    It follows that, this singular kernel satisfies the Assumption $\left(\mathbf{A}^{\nabla \cdot K}\right)$. Furthermore, $K$ satisfies $(\mathbf{A}^{ K})$, see \cite{Simon2}, page 8.
\end{rem}
\begin{rem}
   Assumption \((\mathbf{A}^{\rho_0})\) was formulated in \cite{Feng2,Feng,Feng1}, where, relying on Li–Yau–type gradient bounds and Hamilton–type heat kernel estimates, the authors use it to control the solution and its derivatives up to second order. It allows the relative entropy approach of \cite{Jabin_2018} to be extended from the periodic framework to the whole space. It motivates the following definition.
\end{rem}
\subsection{Solution of limiting process} \label{solution}
We now formulate a suitable notion of solution of (\ref{SPDE_Ito}) to derive our quantitative propagation of chaos. 
\begin{defn} \label{existence sol}
Assume $(\mathbf{A}^{\nabla\cdot K})$,   $(\mathbf{A}^{\sigma})$, and $(\mathbf{A}^{\rho_0})$. A pathwise solution for (\ref{SPDE_Ito}) with initial condition $\rho_0$ is a measurable function $\rho:[0,T]\times \mathbb{R}^d \times \Omega \to \mathbb{R}$, for some $T>0$, such that:
\begin{itemize}
    \item  it holds $\mathbb{P}$-a.s.
    \begin{align}
        \rho \in L^\infty\left([0,T]; W^{2,1}\cap W^{2,\infty} \cap \mathcal{P}(\mathbb{R}^d)\right) \label{regularityW}
    \end{align}
    \item the solution $\rho$ verifies, $\mathbb{P}$-a.s., that
\begin{align}
        \left|\nabla \ln{\rho}(t,x)\right|^2 \leq C_1(1+|x - X_t|^2), \label{decay grad ln rho}
    \end{align}
     \begin{align}
        \left|\nabla^2 \ln{\rho}(t,x)\right| \leq C_2(1+|x - X_t|^2) \label{decay hess ln rho}
    \end{align}
    and
     \begin{align}
       \rho(t,x) \leq \frac{C_3}{1+t}\exp{\left(-\frac{|x - X_t|^2}{8t + C_3}\right)}, \label{decay rho}
    \end{align}
    where
    \begin{align}
        dX_t = \sigma_t dB_t, \label{random}
    \end{align}
    for all $t \in [0,T]$ and constants $C_1, C_2$ and $C_3$ dependent on $\widetilde{C}_1$, $\widetilde{C}_2$ and $\widetilde{C}_3$.
    \end{itemize}
\end{defn}

\subsection{Statement of the main results} \label{mainr}
We are now able to extend what was known before for bounded kernels in the whole space as in \cite{Niko} and in \cite{Ale}, for sub-Coulombian kernels, under periodic boundary conditions, by deriving entropy estimates. 
Also we extend previous results by \cite{Josue} and \cite{Pisa}, by establishing energy estimates. Our
main result is stated as follows.

 \begin{thm} \label{second main}
Assume $(\mathbf{A}^V)$, $(\mathbf{A}^{K})$,  $(\mathbf{A}^{\nabla\cdot K})$,   $(\mathbf{A}^{\sigma})$,  $(\mathbf{A}^{\rho_0})$, $(\mathbf{A}^{X_0})$, $d=2$  and let $\rho$ be a solution of (\ref{SPDE_Ito}), given by Definition  \ref{existence sol}. In addition, let the dynamics of the particle system be given by (\ref{particles})   and  
\begin{align}
\lim_{N\to \infty} N^{\theta} \mathcal{H}(\rho_0^N|\rho_0)  = 0, \, \, \, \mathbb{P}-a.s., \label{initial entropy1}
\end{align}
\begin{align}
\lim_{N\to \infty} N^{\theta }\left\|\rho_0^N - \rho_0\right\|_2^2  = 0, \, \, \, \mathbb{P}-a.s. \label{initial energy}
\end{align}
and
\begin{align*}
    \theta \doteq \min{\left(\left(1 - \beta\left(2 +2\alpha\right)\right); \left(\frac{1}{2}-\beta\Big(1 + \frac{1 }{2}\Big) \right);  \frac{\beta \gamma^2}{2}  \right)} - \delta
\end{align*}
with  $\gamma \in (0,1/2)$ and $\delta > 0$, such that $\theta >0$ and $\beta \in \left(0,\frac{1}{4}\right) \cap \left(0,\frac{1}{\big[2 + 2\alpha\big]}\right)$.
\vspace{.3cm}

Then, there exists a time $T >0$, such that
\begin{align}
    \lim_{N \to \infty}N^{\theta}\left(\sup_{t\in [0,T]}\left(\mathcal{H}(\rho_{t}^N|\rho_{t}) + \left\| \rho_{t}^N - \rho_{t}\right\|_{2}^{2}\right) + \int_0^T \left\|\nabla(\rho_s^N - \rho_s)\right\|_{2}^{2}\, ds \right)&=0, 
     \end{align}
     almost surely.
\end{thm}
As a corollary of Theorem \ref{second main}, we establish quantitative propagation of chaos for the marginals of the empirical measure of the particle system. In the sense of \cite[Section 8.3]{BogachevII}, it corresponds to get convergence of the empirical measure in the sense to the Kantorovich–Rubinstein metric, given in  \eqref{eq:defWasserstein}.		
			
\begin{cor} \label{cor:rateEmpMeas}
Under the assumptions of Theorem~\ref{second main}, we have 
   \begin{align}
  \lim_{N\to \infty} N^{\theta}\sup_{t \in [0,T]}\left\|S_{t}^N - \rho_{t} \right\|_{0}^2 &=0, \, \, \, \mathbb{P}-a.s. \label{grownallffffg}
   \end{align}
\end{cor}
\vspace{.2cm}

In dimensions greater than or equal to two, we establish the boundedness of the relative entropy functional:
  \begin{thm} \label{first main}
Assume $(\mathbf{A}^V)$,  $(\mathbf{A}^{\nabla\cdot K})$,   $(\mathbf{A}^{\sigma})$,  $(\mathbf{A}^{\rho_0})$, $(\mathbf{A}^{X_0})$  and let $\rho$ be a solution of (\ref{SPDE_Ito}), given by Definition  \ref{existence sol}. In addition, let the dynamics of the particle system be given by (\ref{particles})   and  
\begin{align}
\lim_{N\to \infty} N^\theta \mathcal{H}(\rho_0^N|\rho_0)  = 0, \, \, \, \mathbb{P}-a.s. \label{initial entropy}
\end{align}
and
\begin{align*}
\theta = \min \left(1 - \beta\left(1 +\frac{2}{d} + 2\alpha\right);\frac{1}{2}-\beta\Big(1 + \frac{1 }{d}\Big) \right) - \delta 
\end{align*}
with  $\delta > 0$, such that $\theta >0$, $d \geq 1$ and $\beta \in \left(0,\frac{1}{2\big[1 + \frac{1}{d}\big]}\right) \cap \left(0,\frac{1}{\big[1 + \frac{2}{d} + 2\alpha\big]}\right) $.
\vspace{.1cm}

Then, there exist a time $T >0$,  an almost surely finite random variable $A_0$, and  a set $\Lambda \in \mathcal{F}$, $\mathbb{P}(\Lambda)=1$, with the following property: for all $\omega \in \Lambda$, there exists $N_0=N_0(\omega)$ such that, if $N\geq N_0$ we have
\vspace{.1cm}
\begin{align}
   \sup_{t \in [0,T]}\mathcal{H}(\rho_{t}^N|\rho_{t}) &\lesssim  \left(N^{-\theta} + A_0N^{-\theta} + T \right). \label{grownallff}
   \end{align}

   In particular, we get
   \begin{align}
  \limsup_{N\to \infty} \sup_{t \in [0,T]}\mathcal{H}(\rho_{t}^N|\rho_{t}) &\lesssim  T, \, \, \, \, \, \, \mathbb{P}-a.s. \label{grownallffff}
   \end{align}
\end{thm}

 The next result, whose proof is postponed to the Appendix, establishes that the limiting equation is well posed in the appropriate functional framework, given by Definition \ref{existence sol}. 
\begin{thm}\label{existence}
Assume $(\mathbf{A}^{\nabla\cdot K})$,   $(\mathbf{A}^{\sigma})$, and $(\mathbf{A}^{\rho_0})$.
Then, there exists a solution of (\ref{SPDE_Ito}) in the sense of Definition \ref{existence sol}.
\end{thm} 
\section{Proof of Theorem \ref{second main}} \label{Proof}

In this section  we present the proofs of Theorem \ref{second main}.
In Subsection \ref{timeevolution}, by applying the Itô's formula, such as in \cite{Ale}, we derive an evolution equation for the relative entropy functional of the regularized empirical measure with respect to the solution of \eqref{SPDE_Ito}. In
Subsection \ref{quadr}, we deal with the quadratic variation terms arising from the application of Itô’s formula, while
in Subsection \ref{nonlineaa}, we address the nonlinear contributions appearing in the evolution equation for the relative entropy functional obtained in the first step.

In Subsection \ref{timeevolution1}, we derive an evolution equation for the energy, relating our densities, being the nonlinear term, arising from this computation, estimated in Subsection \ref{nonlinea1}, as in Subsection \ref{fin}, we combine our entropy and energy estimates.

Finally, in Subsection  \ref{end of the proof}, by a nonlinear Gronwall Lemma and suitable localization techniques, we close the argument and complete the proof.

\subsection{Time evolution of the relative entropy} 
\label{timeevolution}
By applying Itô’s formula, following the same arguments developed in Section 2 of \cite{Ale}, we derive the corresponding evolution equation for the entropy functional between the regularized empirical measure $\rho^N$ and the solution $\rho$ of (\ref{SPDE_Ito}): 
\begin{align}
    \mathcal{H}(\rho_t^N|\rho_t)- \mathcal{H}(\rho_0^N|\rho_0) 
     &= I_t + II_t + III_t + M_t^N + IV_t,\label{ito at entropy}
\end{align}
where 
\begin{align}
    I_t &\overset{}{=}
     \int_0^t\int_{\mathbb{R}^d} \rho_s^N\nabla \ln\left(\frac{\rho_s}{\rho_s^N} \right)[K\ast \rho_s -  K\ast \rho_s^N] \,dx \,ds  \nonumber\\
     &-
     \int_0^t\int_{\mathbb{R}^d} \nabla \ln\left(\frac{\rho_s}{\rho_s^N} \right) \left\langle S_s^N, V^N(x - \cdot)[K\ast \rho_s^N(\cdot) - K\ast \rho_s^N(x)]\right\rangle \,dx \,ds,  \label{It}
\end{align}
\begin{align}
    II_t + IV_t &= \frac{1}{N^2}\sum_{i=1}^N \int_0^t \int_{\mathbb{R}^d}  \frac{1}{\rho_s^N}|\nabla V^N(x - X_s^{i,N})|^2 \,dx\,ds, \label{quadratic var}
\end{align}
\begin{align}
    III_t &= - \int_0^t\mathcal{I}\left(\rho_s^N|\rho_s\right) \,ds, \label{II_t}
\end{align}
and
\begin{align}
    M_t^N &= \int_0^t \int_{\mathbb{R}^d}\left[\ln(\rho_s) -\ln(\rho_s^N)\right]\sigma^{\top}_s \nabla \rho_s^N  \, dx \, dB_s\nonumber\\
       & + \int_0^t\int_{\mathbb{R}^d}\sigma^{\top}_s \nabla \rho_s^N  \, dx \, dB_s \nonumber\\
     &+ \frac{1}{N} \sum_{i=1}^N \int_0^t \int_{\mathbb{R}^d}  \left[\ln(\rho_s) - \ln(\rho_s^N)\right]\nabla V^N(x - X_s^i) \,dx \,dW^{i,N}_s \nonumber \\
    &- \frac{1}{N} \sum_{i=1}^N  \int_0^t\int_{\mathbb{R}^d} \nabla V^N(x - X_s^i) \,dx \,dW^{i,N}_s  \nonumber \\
     & + \int_0^t \int_{\mathbb{R}^d} \frac{\rho_s^N}{\rho_s} \sigma^{\top}_s\nabla \rho_s  \, dx \, dB_s.  \label{martingalea}
\end{align}

\subsubsection{Estimates for $II_t + IV_t$: quadratic variation terms} \label{quadr}
We now derive estimates for the quadratic variation terms appearing in \eqref{ito at entropy}. The main novelty consists in introducing a suitable stopping time in order to handle the integral over the whole space arising in this term. It is important to emphasize that the argument used in \cite{Ale}, which relies on periodic boundary conditions, is no longer available in the present setting.

First, taking into account Assumption $(\mathbf{A}^V)$ in (\ref{quadratic var}), we have
\begin{align}
     II_t + IV_t 
     &\overset{(\ref{mollifier inequality})}{\leq} C_d\frac{N^{\beta/d}}{N^2}\sum_{i=1}^N \int_0^t\int_{\mathbb{R}^d}  \frac{1}{\rho_s^N}| V^N(x - X_s^{i,N})|^2 \,dx\,ds
     .\label{II + IV}
   \end{align}
  We introduce the following stopping time:
  \begin{align}
      \tau^{N}=\inf \left\{t \geq 0;  \exists i \in \{1,...,N\}, |X_t^{i,N}|\geq N^\beta  \right\} \wedge T. \label{stttoping}
  \end{align}
By decay of $\rho_0$ in (\ref{decayi}), along with Borel-Cantelli Lemma, there exists $\Lambda_0 \in \mathcal{F}$, $\mathbb{P}(\Lambda_0)=1$, with the following property: for all $\omega \in \Lambda_0$, there exists $N_0=N_0(\omega)$ such that $\tau^N(\omega) >0$, for all $N\geq N_0$. 

Indeed, setting
\begin{align}
    \Lambda^N\doteq \left\{\exists i \in \{1,...,N\}, |X_0^{i,N}| \geq N^\beta\right\} = \bigcup_{i=1}^N \left\{ |X_0^{i,N}| \geq N^\beta\right\} \label{set stop}
\end{align}
by Assumption $(\mathbf{A}^{X_0})$, we derive
\begin{align}    \mathbb{P}\left(\Lambda^N\right)&\overset{(\ref{set stop})}{\leq} \sum_{i=1}^N \mathbb{P}\left( |X_0^{i,N}|\geq  N^\beta\right) \nonumber\\
    &= \sum_{i=1}^N \int_{|x|\geq N^\beta}\rho_0(x)\, dx \nonumber\\
    &\leq C_d \sum_{i=1}^N (N^\beta)^{d-2} e^{-\widetilde{C}_3^{-1} (N^\beta)^2} 
    \nonumber\\
    &\leq C_d  N(N^\beta)^{d-2} e^{-\widetilde{C}_3^{-1} (N^\beta)^2} \label{decay initial}
\end{align}
where in the second inequality we are using (\ref{decayi}) 
to get
\begin{align}
   \int_{|x|\ge N^\beta} \rho_0(x) dx &\overset{(\ref{decayi})}{\leq} \widetilde{C}_3\int_{|x|\ge N^\beta} e^{-\widetilde{C}_3^{-1} |x|^2} dx \nonumber\\
    &=C_d \int_{r=N^\beta}^{\infty} e^{-\widetilde{C}_3^{-1}r^2}  r^{d-1}  dr \le C_d (N^\beta)^{d-2} e^{-\widetilde{C}_3^{-1} (N^\beta)^2}.
\end{align}
It follows by Borel-Cantelli Lemma, 
\begin{align}    \mathbb{P}\left(\bigcap_{M=1}^{\infty} \bigcup_{N=M}^\infty \Lambda^N\right) = 0. \label{zero measure}
\end{align}
With this in view, since $t \mapsto X_t^{i,N}$ is continuous, for all $i \in \{1,...,N\}$, almost surely, by setting
\begin{align*}
    \Lambda_{X_0} \doteq \bigcap_{N\in \mathbb{N}}\bigcap_{i=1}^N\left\{t\mapsto X_t^{i,N} \, \, \, \text{is continuous} \right\}
\end{align*}
and
\begin{align}
    \Lambda_0 \doteq \left(\bigcup_{M=1}^{\infty} \bigcap_{N=M}^\infty (\Lambda^N)^c \right) \cap \Lambda_{X_0}\label{medida total}
\end{align}
we have that $\mathbb{P}(\Lambda_0) =1$, and for all $\omega \in \Lambda_0$,
there exists $N_0=N_0(\omega)$, such that $\tau^N(\omega)>0$, for all $N\geq N_0$. 

Let $s \in (0,\tau^N)$. We observe that, 
\begin{align}
    I^N\doteq \int_{\mathbb{R}^d}\frac{|V^N(x-X_s^{i,N})|^2}{\rho^N_s(x)}\, dx &= \int_{|x|\leq 2|X_s^{i,N}|}\frac{|V^N(x-X_s^{i,N})|^2}{\rho^N_s(x)}\, dx \nonumber\\
    &+ \int_{|x|> 2|X_s^{i,N}|}\frac{|V^N(x-X_s^{i,N})|^2}{\rho^N_s(x)}\, dx. \label{split Xi}
\end{align}
Now, if $|x|\leq 2|X_s^{i,N}|$, then $|x|< 2N^\beta$ for $s \in (0,\tau^N)$, and since $V^N\leq N^\beta$, we derive for $\alpha>d/2$
\begin{align}
    V^N(x - X_s^{i,N}) \leq N^\beta \left(1+|x|^2\right)^{\alpha-\alpha} \leq C_{\alpha}  N^{2\alpha \beta + \beta}\left(1+|x|^2\right)^{-\alpha}. \label{split Xi2}
\end{align}
On the other hand, if $|x|> 2|X_s^{i,N}|$ we have $|x - X_s^{i,N}| > \frac{|x|}{2}$ and by decay of $V$, given in (\ref{decay mollifier}), we get 
\begin{align}
    V^N(x - X_s^{i,N}) &\overset{(\ref{decay mollifier})}{\leq} C_\alpha N^{\beta}\left(1+|N^{\frac{\beta}{d}}(x-X_s^{i,N})|^2\right)^{-\alpha}\nonumber\\\
    &\leq C_\alpha N^{\beta}\left(1+\frac{|N^{\frac{\beta}{d}}x|^2}{4}\right)^{-\alpha}\leq C_\alpha N^{\beta}\left(1+\frac{|x|^2}{4}\right)^{-\alpha}. \label{split Xi1}
\end{align}

Therefore by (\ref{split Xi}),  (\ref{split Xi2}) and (\ref{split Xi1}), we find
\begin{align}    
I^N &\leq C_{\alpha} \int_{|x|\leq 2|X_s^{i,N}|}\frac{|V^N(x-X_s^{i,N})|}{\rho^N_s(x)}\left(N^{2\alpha \beta + \beta}\left(1+|x|^2\right)^{-\alpha}\right)\, dx\nonumber\\
    &+ C_{\alpha}\int_{|x|> 2|X_s^{i,N}|}\frac{|V^N(x-X_s^{i,N})|}{\rho^N_s(x)}\left(N^\beta\left(1+\frac{|x|^2}{4}\right)^{-\alpha}\right)\, dx \nonumber\\
    &\leq C_{\alpha}N^{2\alpha\beta + \beta}\int_{\mathbb{R}^d}\frac{|V^N(x-X_s^{i,N})|}{\rho^N_s(x)}\left(1+|x|^2\right)^{-\alpha} \, dx
    \nonumber\\
    &+ C_{\alpha}N^{2\alpha\beta + \beta}\int_{\mathbb{R}^d}\frac{|V^N(x-X_s^{i,N})|}{\rho^N_s(x)} \left(1+\frac{|x|^2}{4}\right)^{-\alpha}\, dx
    \nonumber\\
    &=C_{\alpha} N^{2\alpha \beta + \beta}\int_{\mathbb{R}^d}\frac{|V^N(x-X_s^{i,N})|}{\rho^N_s(x)}\phi(x)\, dx, \label{split Xi4}
\end{align}
with
\begin{align}
    \phi(x) \doteq \left(\left(1+|x|^2\right)^{-\alpha} + \left(1+\frac{|x|^2}{4}\right)^{-\alpha}\right).
\end{align}
Hence by (\ref{split Xi4}) in (\ref{II + IV}) we obtain
\begin{align}
    II_{t \wedge \tau^N} + IV_{t \wedge \tau^N} 
    &\overset{(\ref{II + IV})}{\leq} \frac{N^{\frac{2\beta}{d}}}{2N^2}\int_0^{t \wedge \tau^N}\sum_{i=1}^N \int_{\mathbb{R}^d}  \frac{1}{\rho_s^N(x)}|V^N(x - X_s^{i,N})|^2 \,dx \, ds\nonumber\\
    &\overset{}{\leq}C_{\alpha,d} \frac{N^{\beta(1 +\frac{2}{d} + 2\alpha) }}{N^2}\int_0^{t \wedge \tau^N}\sum_{i=1}^N \int_{\mathbb{R}^d}\frac{|V^N(x-X_s^{i,N})|}{\rho^N_s(x)}\phi(x)\, dx\,ds\nonumber\\
    &\overset{}{\leq}C_{\alpha,d} \frac{N^{\beta(1 +\frac{2}{d} + 2\alpha) }}{N}\int_0^{t \wedge \tau^N} \int_{\mathbb{R}^d}\frac{\frac{1}{N}\sum_{i=1}^N V^N(x-X_s^{i,N})}{\rho^N_s(x)}\phi(x)\, dx\,ds\nonumber\\
    &\overset{}{\leq}C_{\alpha,d} \frac{N^{\beta(1 +\frac{2}{d} + 2\alpha) }}{N}\int_0^{t \wedge \tau^N} \int_{\mathbb{R}^d}\phi(x)\, dx\,ds
    \nonumber\\
    &\overset{\alpha>d/2}{\leq}C_{\alpha,d,T} N^{-\theta_1},\label{cancelations}
\end{align}
 for $\theta_1 \doteq \left(1 - \beta\left(1+2/d+2\alpha\right)\right)$.       
\subsubsection{Estimates for $I_t$: the nonlinearity} \label{nonlineaa}

We now address the nonlinear terms in (\ref{ito at entropy}). The main novelty is to use the stopping time obtained in previous Subsection, along with the  dissipation that comes from of Fisher information and Ladyzenkaya inequality, to tackle with the singularity of kernel $K$.

First, we set
\begin{align}
    I_t &\overset{}{=}
     \int_0^t\int_{\mathbb{R}^d} \rho_s^N\nabla \ln\left(\frac{\rho_s}{\rho_s^N} \right)[K\ast \rho_s -  K\ast \rho_s^N] \,dx \,ds  \nonumber\\
     &-
     \int_0^t\int_{\mathbb{R}^d} \nabla \ln\left(\frac{\rho_s}{\rho_s^N} \right) \left\langle S_s^N, V^N(x - \cdot)[K\ast \rho_s^N(\cdot) - K\ast \rho_s^N(x)]\right\rangle \,dx \,ds \nonumber\\
     &\doteq I_t^1 + I_t^2\label{only Holder}
\end{align}
Thus, by Young inequality together with (\ref{Holderbiot}), we note that
\begin{align}
    I_t^1 &\leq \epsilon \int_0^t \mathcal{I}\left(\rho_s^N|\rho_s\right) \, ds + C_\epsilon \int_0^t\left(\|\rho_s^N - \rho_s\|_1^2 + \|\rho_s^N - \rho_s\|_4^2\right) \, ds\nonumber\\
    &\leq \epsilon \int_0^t\left( \mathcal{I}\left(\rho_s^N|\rho_s\right) + \|\nabla (\rho^N_s - \rho_s)\|_2^2 \right) \, ds \nonumber\\
    &\hspace{80px}+ C_{\epsilon,d} \int_0^t\left(\|\rho_s^N - \rho_s\|_1^2 + \|\rho_s^N - \rho_s\|_2^2\right) \, ds, \label{idd I_1}
\end{align}
where, in the second step, we use the Ladyzhenskaya and Young inequalities.

Now, we will focus on the term $I_t^2$ in (\ref{only Holder}). 
First, we observe that by subtracting and adding $K\ast \rho_s(\cdot)$ and $K\ast \rho_s(x)$ along with (\ref{Holderbiot})
\begin{align}
    \left|K \ast \rho_s^N(\cdot) - K\ast \rho_s^N(x)\right| & \lesssim \left(\|\rho^N_s - \rho_s\|_1 + \|\rho^N_s - \rho_s\|_4 \right) \nonumber\\
    &\hspace{80px}+ \left|K\ast \rho_s(\cdot) - K\ast \rho_s(x)\right|.\label{K_0 Holder}
\end{align}
Additionally, by (\ref{Holderbiot}) and (\ref{decay rho}), we obtain
\begin{align}
     \left\langle S_s^N, V^N(x - \cdot)\left|K\ast \rho_s(\cdot) - K\ast \rho_s(x)\right|\right\rangle  \nonumber\\
     &\hspace{-140px}\overset{(\ref{Holderbiot})}{\lesssim} \left(\|\rho_s\|_1 + \|\rho_s\|_4\right)\left|\left\langle S_s^N, V^N(x - \cdot)|x - \cdot|^\gamma\right\rangle \right|
     \nonumber\\
     &\hspace{-80px}\overset{(\ref{decay rho})}{\lesssim} \left|\left\langle S_s^N, V^N(x - \cdot)|x - \cdot|^\gamma\right\rangle \right|. \label{trian}
\end{align}
We define
\begin{align}
    R \doteq \left\{(x,\cdot); |x - \cdot| < N^{-\beta \gamma/d}\right\} \label{R}
\end{align}
It follows by (\ref{trian}) and (\ref{R}) that
\begin{align}
   \int_{\mathbb{R}^d} \left\langle S_s^N, V^N(x - \cdot)\left|K\ast \rho_s(\cdot) - K\ast \rho_s(x)\right|\mathds{1}_R(x,\cdot)\right\rangle  \lesssim N^{-\beta \gamma^2/d}.\label{ddd}
\end{align}
On the other hand, by decay of $V$, given by (\ref{decay mollifier12}), we have
\begin{align}
    V^N(x - \cdot)|x - \cdot|^\gamma &= N^\beta V(N^{\beta/d}(x-\cdot))|x-\cdot|^\gamma \nonumber\\
    &\overset{(\ref{decay mollifier12})}{\lesssim} N^\beta \left|N^{\beta/d}(x-\cdot)\right|^{-3d/\beta} \left|N^{\beta/d}(x-\cdot)\right|^{-\gamma}|x-\cdot|^\gamma \nonumber\\
    & \lesssim  N^{-3+\beta} \left|(x-\cdot)\right|^{-3d/\beta}.  \label{ddecay}
\end{align}
Hence, we derive
\begin{align}
   \int_{\mathbb{R}^d} \left\langle S_s^N, V^N(x - \cdot)\left|K\ast \rho_s(\cdot) - K\ast \rho_s(x)\right|\mathds{1}_{R^c}(x,\cdot)\right\rangle  \,dx \nonumber\\
   &\hspace{-170px}\overset{(\ref{Holderbiot}), (\ref{decay rho})}{\lesssim}\int_{\mathbb{R}^d} \left|\left\langle S_s^N, V^N(x - \cdot)|x- \cdot|^{\gamma}\mathds{1}_{R^c}(x,\cdot)\right\rangle \right| \,dx \nonumber\\
   &\hspace{-170px}\overset{(\ref{ddecay})}{\lesssim} N^{-3+\beta}\int_{|x| \leq N^\beta} \left|\left\langle S_s^N, \left|(x-\cdot)\right|^{-3d/\beta}\mathds{1}_{R^c}(x,\cdot)\right\rangle \right| \,dx \nonumber\\
   &\hspace{-170px}+N^{-3+\beta}\int_{N^\beta < |x| \leq 2N^\beta} \left|\left\langle S_s^N,\left|(x-\cdot)\right|^{-3d/\beta}\mathds{1}_{R^c}(x,\cdot)\right\rangle \right| \,dx \nonumber\\
   &\hspace{-170px}+N^{-3+\beta}\int_{|x| > 2N^\beta} \left|\left\langle S_s^N, \left|(x-\cdot)\right|^{-3d/\beta}\mathds{1}_{R^c}(x,\cdot)\right\rangle \right| \,dx\nonumber\\
   &\hspace{-90px}\doteq I^1 + I^2 + I^3. \label{dd1}
\end{align}
Now, for $\gamma \in (0,1/2)$, we get
\begin{align}
    I^1 &=  N^{-3+\beta}\int_{|x| \leq N^\beta} \left|\left\langle S_s^N, \left|(x-\cdot)\right|^{-3d/\beta}\mathds{1}_{R^c}(x,\cdot)\right\rangle \right| \,dx \leq N^{-3+\beta} N^{d\beta} N^{3\gamma} \nonumber\\
    &\leq N^{-3 + (1+d)\beta + 3\gamma}. \label{dd2}
\end{align}
Also, for the same reason we arrive at,
\begin{align}
    I^2 &=  N^{-3+\beta}\int_{N^{\beta} < |x| \leq 2 N^\beta} \left|\left\langle S_s^N, \left|(x-\cdot)\right|^{-3d/\beta}\mathds{1}_{R^c}(x,\cdot)\right\rangle \right| \,dx \leq C_d N^{-3+\beta} N^{d\beta} N^{3\gamma} \nonumber\\
    &\leq C_d N^{-3 + (1+d)\beta + 3\gamma}. \label{dd3}
\end{align}
Moreover, if $|x| > 2N^\beta$, since that for $s \in (0,\tau^N)$, $|X_s^{i,N}| < N^\beta$ (see (\ref{set stop})), for all $i \in \{1,...,N\}$,  we have, $|x-\cdot| \geq 1/2|x|$, which implies
\begin{align*}
    \int_{|x|> 2 N^\beta} |x - \cdot|^{-3d/\beta}\, dx \leq C_d \int_{|x|> 2 N^\beta} |x|^{-3d/\beta}\, dx = C_d N^{-3d + d\beta}
\end{align*}
and then,
\begin{align}
  I^3 &= N^{-3+\beta} \int_{\mathbb{R}^d}\int_{|x|> 2 N^\beta} |x - y|^{-3d/\beta}\mathds{1}_{R^c}(x,y)\, dx \, S_s^N(dy)\nonumber\\
  &\leq N^{-3+\beta} C_d (N^{\beta})^{-3d/\beta + d}\leq C_d N^{-3 + (1+d)\beta + 3\gamma}. \label{dd4}
\end{align}
It follows  by (\ref{dd1})-(\ref{dd4}) that,
\begin{align}
   \int_{\mathbb{R}^d} \left\langle S_s^N, V^N(x - \cdot)\left|K\ast \rho_s(\cdot) - K\ast \rho_s(x)\right|\mathds{1}_{R^c}(x,\cdot)\right\rangle  \lesssim N^{-3 + (1+d)\beta + 3\gamma}.\label{ddd1}
\end{align}
\vspace{.3cm}

Therefore, by (\ref{ddd}) and (\ref{ddd1}), $\theta_2 \doteq \beta\gamma^2/d$ with $\gamma \in (0,1/2)$, and $s \in (0,\tau^N)$, we deduce,
\begin{align}
   & \int_{\mathbb{R}^d}\left|\nabla \ln\left(\frac{\rho_s}{\rho_s^N} \right)\right| \left\langle S_s^N, V^N(x - \cdot)\left|K\ast \rho_s^N(\cdot) - K\ast \rho_s^N(x)\right|\right\rangle  \,dx \nonumber\\
    &\overset{(\ref{Holderbiot}), (\ref{decay rho}),(\ref{K_0 Holder})}{\lesssim} \int_{\mathbb{R}^d}\left|\nabla \ln\left(\frac{\rho_s}{\rho_s^N} \right)\right| (\rho_s^N)^{1/2+1/2}\left(\|\rho^N_s - \rho_s\|_1 + \|\rho^N_s - \rho_s\|_4 \right) \,dx \nonumber\\
    &+ \int_{\mathbb{R}^d}\left|\nabla \ln\left(\frac{\rho_s}{\rho_s^N} \right)\right|(\rho_s^N)^{1/2} \left|\left\langle S_s^N, V^N(x - \cdot)[K\ast \rho_s(\cdot) - K\ast \rho_s(x)]\right\rangle \right|^{1/2} \,dx\nonumber\\
    &\overset{(\ref{ddd}), (\ref{ddd1})}{\leq} \epsilon \mathcal{I}\left(\rho_s^N|\rho_s\right) + C_{\epsilon,\gamma,\beta, d} \left(\|\rho^N_s - \rho_s\|_1^2 + \|\rho^N_s - \rho_s\|_4^2 + N^{-\theta_2} \right)
    \nonumber\\
    &\overset{(\ref{ddd}), (\ref{ddd1})}{\leq} \epsilon \left( \mathcal{I}\left(\rho_s^N|\rho_s\right) + \|\nabla (\rho_s^N - \rho_s)\|_2^2 \right) \nonumber\\
    &\hspace{120px}+ C_{\epsilon,\gamma,\beta, d} \left(\|\rho^N_s - \rho_s\|_1^2 + \|\rho^N_s - \rho_s\|_2^2 + N^{-\theta_2} \right), \label{dddI2}
\end{align}
where the last step follows from the Ladyzhenskaya and Young inequalities.

Finally, by (\ref{dddI2}), we get
\begin{align}
    I_{t\wedge \tau^N}^2 &\leq 
    \int_0^{t\wedge \tau^N} \epsilon \left( \mathcal{I}\left(\rho_s^N|\rho_s\right) + \|\nabla (\rho_s^N - \rho_s)\|_2^2 \right) \, ds \nonumber\\
    &\hspace{80px}+ C_{\epsilon,\gamma,\beta, d} \int_0^{t\wedge \tau^N}\left(\|\rho^N_s - \rho_s\|_1^2 + \|\rho^N_s - \rho_s\|_2^2 + N^{-\theta_2} \right) \, ds, \label{idd I_2}
\end{align}
and thus, by (\ref{only Holder}), (\ref{idd I_1}) and (\ref{idd I_2}), we end up with
\begin{align}
    I_{t\wedge \tau^N} &\leq 
    2 \epsilon\int_0^{t\wedge \tau^N}  \left( \mathcal{I}\left(\rho_s^N|\rho_s\right) + \|\nabla (\rho_s^N - \rho_s)\|_2^2 \right) \, ds \nonumber\\
    &\hspace{80px}+ C_{\epsilon,\gamma,\beta, d} \int_0^{t\wedge \tau^N}\left(\|\rho^N_s - \rho_s\|_1^2 + \|\rho^N_s - \rho_s\|_2^2 + N^{-\theta_2} \right) \, ds, \label{1I1}
\end{align}
where $\theta_2 \doteq \beta\gamma^2/d$ and $\gamma \in (0,1/2)$.
\subsection{Time evolution of energy} \label{timeevolution1}

Applying Itô's lemma, following the approach of \cite{Josue}, we obtain
	\begin{align}
		\left\| \rho_t^N - \rho_t\right\|_{2}^{2}&\le \left\| \rho_0^N -\rho_0\right\|_{2}^{2}+\widetilde{I}^1_t  -\widetilde{M}_{t}^N -2  \int_0^t\left\| \nabla (\rho_s^N - \rho_s)\right\|_2^{2}\, ds  \nonumber\\
        &+T N^{-\widetilde{\theta}_1}  \, \|\nabla V\|_{2}^2,    \label{dec7second}
	\end{align}
	\noindent with $\widetilde{\theta}_1 \doteq \left(1 - 2\beta(1 + \frac{1 }{d})+\beta\right)$, 
	\begin{align*}
		\widetilde{I}^1_t \doteq 2\int_{0}^{t}\int_{\mathbb{R}^d}\nabla (\rho_s& -\rho_s^N)\Big[\rho_sK\ast \rho_s -\left\langle S_{s}^{N}, V^{N}(x-\cdot) K\ast \rho^N_s(\cdot)\right\rangle \Big]\, dx\,ds,
	\end{align*}
	and
    \begin{align*}
        \widetilde{M}^N_{t} \doteq \frac{2\sqrt{2}}{N} \sum_{i=1}^{N}\int_{0}^{t} \int_{\mathbb{R}^d} (\rho_s- \rho_s^N) \,  ( \nabla V^{N})(x-X_{s}^{i,N}) \,dx \, dW_{s}^{i,N}.
    \end{align*}
    \subsubsection{Estimates for $\widetilde{I}_t$: the nonlinearity} \label{nonlinea1}
We now address the nonlinearity in the evolution equation of energy. The main novelty is identifying, within the energy estimates, the Fisher information of the particle system alongside Ladyzhenskaya inequality, to address the nonlinearity.
    
By subtracting and adding the term $\rho_s^N K\ast \rho_s^N$, we have 
\begin{align}
   &\widetilde{I}^1_t = - 2 \int_{\mathbb{R}^2} \nabla (\rho_s -\rho_s^N)\left[\rho_s  K\ast \rho_s - \rho_s^N K\ast \rho^N_s\right]\, dx
    \nonumber\\ 
    & + 2 \int_{\mathbb{R}^2} \nabla (\rho_s -\rho_s^N)\left\langle S_{s}^{N}, V^{N}(x-\cdot)\left[K\ast \rho_s^{N}(\cdot) - K\ast \rho_s^N(x)\right]\right\rangle \, dx \nonumber\\
    &\doteq \widetilde{I}^{1,1} + \widetilde{I}^{1,2}. \label{energypre}
    \end{align}
    In addition, since $\nabla \cdot K = 0$, by subtracting and adding the term $\rho_s$, we find
    \begin{align}
        \left|\widetilde{I}^{1,1} \right|&\leq \left| \underbrace{2\int_{\mathbb{R}^2} \nabla (\rho_s - \rho_s^N)(\rho_s - \rho_s^N) K\ast \rho_s\, dx}_{=0}\right| \nonumber\\
        & + \left| \underbrace{2\int_{\mathbb{R}^2} \nabla (\rho_s - \rho_s^N)(\rho_s - \rho_s^N) K\ast (\rho_s - \rho_s^N)\, dx}_{=0}\right| \nonumber\\
        & + \left| 2 \int_{\mathbb{R}^2} \nabla (\rho_s - \rho_s^N)\rho_s K\ast (\rho_s - \rho_s^N)\, dx \right| \nonumber\\
        &\overset{\text{Holder}+ (\ref{Holderbiot})}{\lesssim}\left(\|\rho_s - \rho_s^N\|_1 + \|\rho_s - \rho_s^N\|_4\right)\|\nabla(\rho_s - \rho_s^N)\|_2\|\rho_s\|_2. \label{energy}
    \end{align}
    By Young inequality, we observe that
    
    \begin{align}
        \|\rho_s - \rho_s^N\|_1 \|\nabla(\rho_s - \rho_s^N)\|_2\|\rho_s\|_2 \leq C_\epsilon\|\rho_s\|_2^2\|\rho_s - \rho_s^N\|_1^2 + \epsilon\|\nabla(\rho_s -\rho_s^N)\|_2^2. \label{energy1}
    \end{align}
    By Ladyzhenskaya and Young inequalities, we deduce
    
     \begin{align}
        &\|\rho_s - \rho_s^N\|_4 \|\nabla(\rho_s - \rho_s^N)\|_2\|\rho_s\|_2 \nonumber\\
        &\hspace{80px}\lesssim \|\rho_s - \rho_s^N\|_2^{1/2}\|\nabla(\rho_s - \rho_s^N)\|_2^{1/2} \|\nabla(\rho_s - \rho_s^N)\|_2\|\rho_s\|_2\nonumber\\
        &\hspace{80px}\leq  C_\epsilon\|\rho_s\|_2^4\|\rho_s - \rho_s^N\|_2^2 + \epsilon\|\nabla(\rho_s -\rho_s^N)\|_2^2. \label{energy2}
    \end{align}
     Therefore by  (\ref{energy1}) and (\ref{energy2}) in (\ref{energy}), we get
     \begin{align}
         \widetilde{I}^{1,1} \leq C_\epsilon\|\rho_s\|_2^2\|\rho_s - \rho_s^N\|_1^2 + C_\epsilon\|\rho_s\|_2^4\|\rho_s - \rho_s^N\|_2^2 + 2\epsilon\|\nabla(\rho_s -\rho_s^N)\|_2^2. \label{energy I11}
     \end{align}
\vspace{.1cm}

\noindent Now we will take care about $\widetilde{I}^{1,2}$ in (\ref{energypre}). 

First, by decay of $V$, given by (\ref{decay mollifier12}), we have, for $|x-\cdot|\geq N^{-\beta \gamma/d}$
\begin{align}
    V^N(x - \cdot)|x - \cdot|^\gamma &= N^\beta V(N^{\beta/d}(x-\cdot))|x-\cdot|^\gamma \nonumber\\
    &\overset{(\ref{decay mollifier12})}{\lesssim} N^\beta \left|N^{\beta/d}(x-\cdot)\right|^{-3d/\beta} \left|N^{\beta/d}(x-\cdot)\right|^{-\gamma}|x-\cdot|^\gamma \nonumber\\
    & \lesssim  N^{-3+\beta} \left|(x-\cdot)\right|^{-3d/\beta}\nonumber\\
    &\lesssim N^{-\widetilde{\theta}_2},  \label{ddecay9}
\end{align}
which implies that
     \begin{align}
      \left|\left\langle S^N_s, V^N(x - \cdot)|x-\cdot|^{\gamma} \right \rangle \right|\leq C_dN^{-\widetilde{\theta}_2} + N^{-\widetilde{\theta}_3}\rho_s^N(x) \label{triangularization  without compact support}
  \end{align}
  with $\widetilde{\theta}_2 \doteq 3-\beta - 3\gamma$ and $\widetilde{\theta}_3 \doteq \frac{\beta \gamma^2}{d}$, for $\gamma \in (0,1/2)$.
\vspace{.1cm}

It follows by (\ref{Holderbiot}), that
\begin{align}    \left|\widetilde{I}^{1,2}\right| &\overset{(\ref{Holderbiot})}{\lesssim} 2 (\|\rho_s^N\|_1 + \|\rho_s^N\|_4) \int_{\mathbb{R}^2}|\nabla(\rho_s - \rho_s^N)|\left\langle S^N_s, V^N(x - \cdot)|x-\cdot|^{\gamma} \right \rangle \, dx \nonumber\\
 &\overset{(\ref{triangularization  without compact support})}{\lesssim} 2 (\|\rho_s^N\|_1 + \|\rho_s^N\|_4) \int_{\mathbb{R}^2}|\nabla(\rho_s - \rho_s^N)|\left(N^{-\widetilde{\theta}_2} + N^{-\widetilde{\theta}_3}\rho_s^N\right) \, dx \nonumber\\
 &\lesssim 2 \|\rho_s^N\|_1  \int_{\mathbb{R}^2}|\nabla(\rho_s - \rho_s^N)|\left(N^{-\widetilde{\theta}_2} + N^{-\widetilde{\theta}_3}\rho_s^N\right) \, dx \nonumber\\
 &+ 2  \|\rho_s^N\|_4 \int_{\mathbb{R}^2}|\nabla(\rho_s - \rho_s^N)|\left(N^{-\widetilde{\theta}_2} + N^{-\widetilde{\theta}_3}\rho_s^N\right)\,dx \nonumber\\
 &\doteq \widetilde{I}^{1,2,1} + \widetilde{I}^{1,2,2} \label{I01}
\end{align}

Now, by subtracting and adding $\rho_s$ along side Young inequality, we note that
\begin{align}
    \widetilde{I}^{1,2,1} &\leq N^{-\widetilde{\theta}_2}\left(\|\nabla \rho_s\|_1 + \|\nabla \rho_s^N\|_1\right) \nonumber\\
    &+ N^{-\widetilde{\theta}_3} \int_{\mathbb{R}^2}|\nabla (\rho_s - \rho_s^N)||\rho_s - \rho_s^N|\, dx
    \nonumber\\
    &+ N^{-\widetilde{\theta}_3} \int_{\mathbb{R}^2}|\nabla (\rho_s - \rho_s^N)||\rho_s|\, dx\nonumber\\
    &\overset{\text{Young}}{\leq} N^{-\widetilde{\theta}_2}\left(\|\nabla \rho_s\|_1 + \|\nabla \rho_s^N\|_1\right) \nonumber\\
    &+ N^{-\widetilde{\theta}_3}\left( \epsilon \|\nabla (\rho_s - \rho_s^N)\|_2^2 + C_\epsilon \|\rho_s - \rho_s^N\|_2^2\right)
    \nonumber\\
    &+ N^{-\widetilde{\theta}_3}\left( \epsilon\|\nabla (\rho_s - \rho_s^N)\|_2^2 + C_\epsilon\|\rho_s\|_2^2\right). \label{I02}
\end{align}

Again, by by subtracting and adding $\rho_s$, we find
\begin{align}
    \widetilde{I}^{1,2,2}& \leq N^{-\widetilde{\theta}_2}\|\rho_s^N\|_4\|\nabla(\rho_s - \rho_s^N)\|_1 + N^{-\widetilde{\theta}_3}\|\rho_s^N\|_4\int_{\mathbb{R}^2}|\nabla(\rho_s -\rho_s^N)| \rho_s^N \, dx \nonumber\\
    & \overset{\text{Holder}}{\leq} N^{-\widetilde{\theta}_2}\|\rho_s\|_4\|\nabla(\rho_s - \rho_s^N)\|_1 + N^{-\widetilde{\theta}_2}\|\rho_s^N - \rho_s\|_4\|\nabla(\rho_s - \rho_s^N)\|_1 \nonumber\\
    &+ N^{-\widetilde{\theta}_3}\|\rho_s - \rho_s^N\|_4\|\nabla(\rho_s - \rho_s^N)\|_2\|\rho_s - \rho_s^N\|_2
    \nonumber\\
    &+ N^{-\widetilde{\theta}_3}\|\rho_s\|_4\|\nabla(\rho_s - \rho_s^N)\|_2\|\rho_s - \rho_s^N\|_2 \nonumber\\
    &+N^{-\widetilde{\theta}_3}\|\rho_s - \rho_s^N\|_4\|\nabla(\rho_s - \rho_s^N)\|_2\|\rho_s\|_2 \nonumber\\
    &+N^{-\widetilde{\theta}_3}\|\rho_s\|_4\|\nabla(\rho_s - \rho_s^N)\|_2\|\rho_s\|_2 \nonumber\\
    &\doteq \widetilde{I}^{1,2,2}_1 + \widetilde{I}^{1,2,2}_2 + \widetilde{I}^{1,2,2}_3 + \widetilde{I}^{1,2,2}_4 + \widetilde{I}^{1,2,2}_5. \label{I03}
\end{align}
Now by Ladyzhenskaya and Young inequalities, we find
\begin{align}
    \widetilde{I}^{1,2,2}_1 &\overset{\text{Young}}{\leq} N^{-\widetilde{\theta}_2}\left(\|\rho\|_4^2 + \|\nabla(\rho_s - \rho_s^N)\|_1^2 + \|\rho_s - \rho_s^N\|_4^2 \right) \nonumber\\
    &\overset{\text{Young}+ \text{Lady}}{\leq} N^{-\widetilde{\theta}_2}\left(\|\rho_s\|_4^2 + \|\nabla(\rho_s - \rho_s^N)\|_1^2 \right) \nonumber\\
    &+ N^{-\widetilde{\theta}_2}\left(C_\epsilon\|\rho_s - \rho_s^N\|_2^2 +\epsilon\|\nabla(\rho_s - \rho_s^N)\|_2^2 \right). \label{I04}
\end{align}
Analogously we deduce
\begin{align}
    \widetilde{I}^{1,2,2}_3 + \widetilde{I}^{1,2,2}_4 + \widetilde{I}^{1,2,2}_5 &\leq N^{-\widetilde{\theta}_3}\left(\epsilon\|\nabla(\rho_s - \rho_s^N\|_2^2 + C_\epsilon\|\rho_s\|_4^2\|\rho_s - \rho_s^N\|_2^2\right)\nonumber\\
    &+ N^{-\widetilde{\theta}_3}\left(\epsilon\|\nabla(\rho_s - \rho_s^N\|_2^2 + C_\epsilon\|\rho_s\|_2^4\|\rho_s - \rho_s^N\|_2^2\right)
    \nonumber\\
    &+ N^{-\widetilde{\theta}_3}\left(\epsilon\|\nabla(\rho_s - \rho_s^N\|_2^2 + C_\epsilon\|\rho_s\|_4^4 + C_\epsilon\|\rho_s\|_2^4\right). \label{I05}
\end{align}
Also, for the term $\widetilde{I}^{1,2,2}_2$, by means of Ladyzhenskaya and Young inequalities we get
\begin{align}
    \widetilde{I}^{1,2,2}_2 \leq N^{-\widetilde{\theta}_3}\left(C_\epsilon\|\rho_s - \rho_s^N\|_2^6 + \epsilon\|\nabla(\rho_s - \rho_s^N)\|_2^2\right). \label{I06}
\end{align}
It follows by (\ref{I02}), (\ref{I03}), (\ref{I04}), (\ref{I05}) and (\ref{I06}) in (\ref{I01}), we arrive at
\begin{align}
    \widetilde{I}^{1,2} 
    &\overset{}{\leq} N^{-\widetilde{\theta}_2}\left(\|\nabla \rho_s\|_1 + \|\nabla \rho_s^N\|_1\right) \nonumber\\
    &+ N^{-\widetilde{\theta}_3}\left( \epsilon \|\nabla (\rho_s - \rho_s^N)\|_2^2 + C_\epsilon \|\rho_s - \rho_s^N\|_2^2\right)
    \nonumber\\
    &+ N^{-\widetilde{\theta}_3}\left( \epsilon\|\nabla (\rho_s - \rho_s^N)\|_2^2 + C_\epsilon\|\rho_s\|_2^2\right) \nonumber\\
    &\overset{}{+} N^{-\widetilde{\theta}_2}\left(\|\rho\|_4^2 + \|\nabla(\rho_s - \rho_s^N)\|_1^2 \right) \nonumber\\
    &+ N^{-\widetilde{\theta}_2}\left(C_\epsilon\|\rho_s - \rho_s^N\|_2^2 +\epsilon\|\nabla(\rho_s - \rho_s^N)\|_2^2 \right) \nonumber\\
    &+ N^{-\widetilde{\theta}_3}\left(\epsilon\|\nabla(\rho_s - \rho_s^N\|_2^2 + C_\epsilon\|\rho_s\|_4^2\|\rho_s - \rho_s^N\|_2^2\right)\nonumber\\
    &+ N^{-\widetilde{\theta}_3}\left(\epsilon\|\nabla(\rho_s - \rho_s^N\|_2^2 + C_\epsilon\|\rho_s\|_2^4\|\rho_s - \rho_s^N\|_2^2\right)
    \nonumber\\
    &+ N^{-\widetilde{\theta}_3}\left(\epsilon\|\nabla(\rho_s - \rho_s^N\|_2^2 + C_\epsilon\|\rho_s\|_4^4 + C_\epsilon\|\rho_s\|_2^4\right) \nonumber\\
    &+N^{-\widetilde{\theta}_3}\left(C_\epsilon\|\rho_s - \rho_s^N\|_2^6 + \epsilon\|\nabla(\rho_s - \rho_s^N)\|_2^2\right). \label{by}
\end{align}

Hence, by (\ref{energy I11}) and (\ref{by}) in (\ref{energypre}), we conclude with
\begin{align}
    \widetilde{I}_t &\leq C_\epsilon\|\rho_s\|_2^2\|\rho_s - \rho_s^N\|_1^2 + C_\epsilon\|\rho_s\|_2^4\|\rho_s - \rho_s^N\|_2^2 + 2\epsilon\|\nabla(\rho_s -\rho_s^N)\|_2^2 \nonumber\\
    &\overset{}{+} N^{-\widetilde{\theta}_2}\left(\|\nabla \rho_s\|_1 + \|\nabla \rho_s^N\|_1\right) \nonumber\\
    &+ N^{-\widetilde{\theta}_3}\left( \epsilon \|\nabla (\rho_s - \rho_s^N)\|_2^2 + C_\epsilon \|\rho_s - \rho_s^N\|_2^2\right)
    \nonumber\\
    &+ N^{-\widetilde{\theta}_3}\left( \epsilon\|\nabla (\rho_s - \rho_s^N)\|_2^2 + C_\epsilon\|\rho_s\|_2^2\right) \nonumber\\
    &\overset{}{+} N^{-\widetilde{\theta}_2}\left(\|\rho\|_4^2 + \|\nabla(\rho_s - \rho_s^N)\|_1^2 \right) \nonumber\\
    &+ N^{-\widetilde{\theta}_2}\left(C_\epsilon\|\rho_s - \rho_s^N\|_2^2 +\epsilon\|\nabla(\rho_s - \rho_s^N)\|_2^2 \right) \nonumber\\
    &+ N^{-\widetilde{\theta}_3}\left(\epsilon\|\nabla(\rho_s - \rho_s^N\|_2^2 + C_\epsilon\|\rho_s\|_4^2\|\rho_s - \rho_s^N\|_2^2\right)\nonumber\\
    &+ N^{-\widetilde{\theta}_3}\left(\epsilon\|\nabla(\rho_s - \rho_s^N\|_2^2 + C_\epsilon\|\rho_s\|_2^4\|\rho_s - \rho_s^N\|_2^2\right)
    \nonumber\\
    &+ N^{-\widetilde{\theta}_3}\left(\epsilon\|\nabla(\rho_s - \rho_s^N\|_2^2 + C_\epsilon\|\rho_s\|_4^4 + C_\epsilon\|\rho_s\|_2^4\right) \nonumber\\
    &+N^{-\widetilde{\theta}_3}\left(C_\epsilon\|\rho_s - \rho_s^N\|_2^6 + \epsilon\|\nabla(\rho_s - \rho_s^N)\|_2^2\right). \label{tilde I}
\end{align}
\subsection{Combining entropy and energy} \label{fin}
First, with the same computations as in \cite{Josue}, for all $\delta \in (0,\theta_3)$ there exists a random variable \( A_0 \),  finite almost surely, such that
\begin{align}       \sup_{t\in[0,T]}|M_t^N|  \leq A_0N^{-\theta_3 + \delta}, \label{marga}
    \end{align} 
with $\theta_3 \doteq \left(\frac{1}{2}-\beta\Big(1 + \frac{1 }{d}\Big) \right)$.

In addition, with the same computations as in \cite{Ale}, for all $\delta \in (0,\widetilde{\theta}_4)$ there exists a random variable \( \widetilde{A}_0 \),  finite almost surely, such that
\begin{align}       \sup_{t\in[0,T]}|\widetilde{M}_t^N|  \leq \widetilde{A}_0N^{-\widetilde{\theta}_4 + \delta}, \label{marga1}
    \end{align} 
with $\widetilde{\theta}_4 \doteq \left(\frac{1}{2}-\beta\Big(1 + \frac{1 }{d}\Big) \right)$.

Now we put $(\ref{II_t}), (\ref{cancelations})$,  (\ref{1I1}) and (\ref{marga}) into $(\ref{ito at entropy})$, and use by Csiszár-Kullback-Pinsker inequality (\ref{pinsker}), to get
\begin{align}
    \mathcal{H}(\rho_{t\wedge \tau^N}^N|\rho_{t\wedge \tau^N}) - \mathcal{H}(\rho_0^N|\rho_0) &\leq C\int_0^{t\wedge \tau^N}\left(\mathcal{H}(\rho_s^N|\rho_s) + \|\rho^N_s - \rho_s\|_2^2\right)\, ds \nonumber\\
    &\overset{}{+}(-1+ 2\epsilon)\int_0^{t\wedge \tau^N} \left(\mathcal{I}(\rho_s^N|\rho_s) + \|\nabla(\rho_s^N - \rho_s)\|_2^2\right)\,ds\nonumber\\    
    &+C N^{-\theta_1} + C N^{-\theta_2} + A_0N^{-\theta_3 + \delta},
    \label{entropy regular grownall}
\end{align}
for a universal constant $C >0$.

Analogously, by (\ref{tilde I}) and (\ref{marga1})   in (\ref{dec7second}), we obtain
\begin{align}
		\left\| \rho_t^N - \rho_t\right\|_{2}^{2} - \left\| \rho_0^N -\rho_0\right\|_{2}^{2} 
        &\leq C_\epsilon \int_0^t\|\rho_s\|_2^2\mathcal{H}(\rho_s^N|\rho_s)  \, ds \nonumber\\
        &+ \int_0^t\left(C_\epsilon + \|\rho_s\|_4^2 + \|\rho_s\|_2^4  \right)\|\rho_s - \rho_s^N\|_2^2 \, ds\nonumber\\
    &+C_\epsilon \int_0^t\|\rho_s - \rho_s^N\|_2^6  \, ds
        \nonumber\\
    &\overset{}{+} 2N^{-\widetilde{\theta}_2}\int_0^t\left(\|\nabla \rho_s\|_1^2 + \|\nabla \rho_s^N\|_1^2 + T\right) \, ds\nonumber\\
    &+ N^{-\widetilde{\theta}_3}  C_\epsilon\int_0^t\|\rho_s\|_2^2 \, ds + N^{-\widetilde{\theta}_2}\int_0^t\|\rho_s\|_4^2 \, ds\nonumber\\
     &+ N^{-\widetilde{\theta}_3}\int_0^t\left( C_\epsilon\|\rho_s\|_4^4 + C_\epsilon\|\rho_s\|_2^4\right) \, ds \nonumber\\
     &+ (-2 +8\epsilon) \int_0^t \left\|\nabla(\rho_s^N - \rho_s)\right\|_{2}^{2}\, ds\nonumber\\
    &+ TN^{-\widetilde{\theta}_1}  \, \|\nabla V\|_{2}^2 \nonumber\\
    &+\widetilde{A}_0N^{-\widetilde{\theta}_4 + \delta}.
    \label{dec7second1}
	\end{align}
    Now, since  $\rho_s^N \in \mathcal{P}\left(\mathbb{R}^d\right)$, $\mathbb{P}$-a.s. and then by Holder's inequality we obtain 
\begin{align}
   \int_{\mathbb{R}^d}|\nabla \rho_s^N| \,dx &= \int_{\mathbb{R}^d}|\nabla \rho_s^N| \frac{\sqrt{\rho_s^N}}{\sqrt{\rho_s^N}} \,dx \nonumber\\
    &\overset{\text{Holder}}{\leq} \left(\int_{\mathbb{R}^d}\frac{|\nabla \rho_s^N|^2}{\rho_s^N} \,dx\right)^{\frac{1}{2}} \left(\int_{\mathbb{R}^d}\rho_s^N \,dx \right)^{\frac{1}{2}} \nonumber \\
    &\overset{}{=} \left(\int_{\mathbb{R}^d}\frac{|\nabla \rho_s^N|^2}{\rho_s^N}\,dx \right)^{\frac{1}{2}}. \label{fisher and L1}
\end{align}
Also, from  Leibniz's rule we deduce, 
\begin{align} \label{fisher rho^N and difference}
   \int_{\mathbb{R}^d}  \frac{|\nabla \rho_s^N|^2}{\rho_s^N} \,dx \nonumber& = \int_{\mathbb{R}^d}  \frac{\left|\nabla \left(\rho_s\frac{\rho_s^N}{\rho_s}\right)\right|^2}{\rho_s^N} \,dx\\
    \nonumber&\leq2\int_{\mathbb{R}^d}  \frac{\left|\nabla \rho_s\left(\frac{\rho_s^N}{\rho_s}\right)\right|^2}{\rho_s^N} \,dx + 2\int_{\mathbb{R}^d}  \frac{\left|\rho_s\nabla \left(\frac{\rho_s^N}{\rho_s}\right)\right|^2}{\rho_s^N} \,dx
    \\
   &\overset{}{\leq} 2\int_{\mathbb{R}^d} \rho_s^N  \frac{\left|\nabla \rho_s\right|^2}{\rho_s^2} \,dx + 2\int_{\mathbb{R}^d} \rho_s^N \left|\nabla \ln \left(\frac{\rho_s^N}{\rho_s}\right)\right|^2 \,dx
    \nonumber\\
      &\overset{}{\leq} 2\eta \mathcal{H}(\rho_s^N|\rho_s)  + C_\eta + 2\int_{\mathbb{R}^d} \rho_s^N \left|\nabla \ln \left(\frac{\rho_s^N}{\rho_s}\right)\right|^2 \,dx, 
   \end{align}
where in the last step we are using the following:
for $\eta > 2C_1(8T + C_3) + 1$
    \begin{align*}        \int_{\mathbb{R}^d}\frac{C_3}{1+t}\exp{\left(-\frac{|x - X_t|^2}{8t + C_3}\right)}\exp{\left(\frac{C_1(1 + |x - X_t|^2)}{\eta}\right)}\, dx < \infty 
    \end{align*}
     which implies by Donsker-Varadhan inequality (\ref{donsker-vara}),
    \begin{align}
    \int_{\mathbb{R}^d} \rho_s^N  \frac{\left|\nabla \rho_s\right|^2}{\rho_s^2} \,dx &\overset{(\ref{donsker-vara})}{\leq} \eta \mathcal{H}(\rho_s^N|\rho_s) + \eta \ln{\left(\int_{\mathbb{R}^d}\rho_s\exp{\left(\frac{|\nabla \ln{\rho_s}|^2}{\eta}\right)}\, dx\right)} \nonumber \\
    &\overset{(\ref{decay grad ln rho}),(\ref{decay rho}), (\ref{logfinite})}{\leq} \eta \mathcal{H}(\rho_s^N|\rho_s) + C_\eta. \label{grad donsker}
\end{align}
Hence since  (\ref{decay rho})  implies $\|\rho\|_{L^\infty L^2}, \|\rho\|_{L^\infty L^4} \leq C_0$, almost surely, by (\ref{dec7second1}), (\ref{fisher and L1})  and (\ref{grad donsker}), there exists a universal constant $C>0$ in (\ref{dec7second1}), we derive
\begin{align}
		&\left\| \rho_t^N - \rho_t\right\|_{2}^{2} +  \int_0^t \left\|\nabla(\rho_s^N - \rho_s)\right\|_{2}^{2}\, ds \nonumber\\
        &\hspace{100px}\leq \left\| \rho_0^N -\rho_0\right\|_{2}^{2} \nonumber\\
        &\hspace{100px}+  \int_0^t C\mathcal{H}(\rho_s^N|\rho_s)  \, ds + \int_0^t C\|\rho_s - \rho_s^N\|_2^2 \, ds\nonumber\\
        &\hspace{100px}+ \int_0^tC\|\rho_s - \rho_s^N\|_2^6  \, ds
        \nonumber\\
   &\hspace{100px}\overset{}{+} CN^{-\widetilde{\theta}_2}\int_0^t\|\nabla \rho_s\|_1^2  \, ds + TCN^{-\widetilde{\theta}_2} \nonumber\\
      &\hspace{100px}\overset{}{+} CN^{-\widetilde{\theta}_2}\int_0^t \left(2\eta \mathcal{H}(\rho_s^N|\rho_s)  + 1 + 2\mathcal{I}(\rho_s^N|\rho_s) \right)\, ds \nonumber\\
      &\hspace{100px}\overset{}{+} CN^{-\widetilde{\theta}_1} + CN^{-\widetilde{\theta}_2}   + CN^{-\widetilde{\theta}_3}   \nonumber\\
    &\hspace{100px}+ CN^{-\widetilde{\theta}_1}  \, \|\nabla V\|_{2}^2 \nonumber\\    &\hspace{100px}+C\widetilde{A}_0N^{-\widetilde{\theta}_4 + \delta}.
    \label{dec7second122}
	\end{align}
It follows that, by adding (\ref{entropy regular grownall}) and (\ref{dec7second122}), for $\epsilon <<1$ and  $N>>1$, we deduce
\begin{align}
    &\mathcal{H}(\rho_{t\wedge \tau^N}^N|\rho_{t\wedge \tau^N}) + \left\| \rho_{t\wedge\tau^N} - \rho^N_{t\wedge\tau^N}\right\|_{2}^{2} + \int_0^{t\wedge\tau^N} \left\|\nabla(\rho_s - \rho^N_s)\right\|_{2}^{2}\, ds \nonumber\\
    &\leq  \mathcal{H}(\rho_0^N|\rho_0) + \|\rho_0 -\rho_0^N\|_2^2 + C_N  \nonumber\\
    &+\int_0^{t\wedge \tau^N}C\left(\mathcal{H} (\rho_s^N|\rho_s) + \|\rho_s - \rho_s^N\|_2^2 +  \left(\int_0^s \left\|\nabla(\rho_r - \rho^N_r)\right\|_{2}^{2}\, dr\right)\right)\, ds \nonumber\\
    & + \int_0^{t\wedge \tau^N}C\left(\mathcal{H}(\rho_s^N|\rho_s) + \|\rho_s - \rho_s^N\|_2^2 + \left(\int_0^s \left\|\nabla(\rho_r - \rho^N_r)\right\|_{2}^{2}\, dr\right)\right)^3  \, ds, \label{quasienergygro}
     \end{align}
     with
     \begin{align*}
      &C_N \doteq   CN^{\widetilde{\theta}}\|\nabla \rho\|_{L^{\infty}L^1}^2  + CN^{\widetilde{\theta}} + CA_0 N^{\widetilde{\theta}} + C\widetilde{A}_0 N^{\widetilde{\theta}}, 
    \end{align*}
and 
\begin{align*}
    \widetilde{\theta} \doteq \min{\left(\theta_1,\theta_2, \theta_3, \widetilde{\theta}_1,\widetilde{\theta}_2, \widetilde{\theta}_3, \widetilde{\theta}_4\right)} - \delta
\end{align*}
for $\delta <<1$.

\subsection{End of the proof: Removal the stopping time}\label{end of the proof}

We now remove the stopping time $\tau^N$ introduced in (\ref{stttoping}) by exploiting the probabilistic framework provided by Assumptions $\left(\mathbf{A}^{\rho_0}\right)$ and $\left(\mathbf{A}^{X_0}\right)$. It allows us to pass from the localized estimates to uniform bounds with high probability, removing the stopping time.

First, since $K=\nabla \cdot K_0$, $K_0 \in L^\infty$ and $V^N = N^\beta V(N^{\beta/d}\cdot)$, by (\ref{particles}), we have
\begin{align}
    \sup_{t \in [0,T]}|X_t^{i,N}| &\leq |X_0^{i,N}| + T\|\nabla V\|_1\|K_0\|_\infty N^\beta \nonumber\\
    &+ \sup_{t\in[0,T]}\left|\int_0^t \sqrt{2}dW_s^{i,N}\right| +  \sup_{t\in[0,T]}\left|\int_0^t \sigma_sdB_s\right|, 
\end{align}
which implies for $\widetilde{C}_{K_0,V} \doteq (1-T\|\nabla V\|_1\|K_0\|_\infty)$
\begin{align}
    \{\tau^N < T\} &\subset \bigcup_{i=1}^N \left\{\sup_{t\in [0,T]} |X_t^{i,N}| \geq N^\beta\right\} \nonumber\\
    &= \left\{\exists i ;\sup_{t\in [0,T]} |X_t^{i,N}| \geq N^\beta\right\}
    \nonumber\\
    &\subset \left\{\exists i;|X_0^{i,N}| + \sup_{t\in [0,T]} \left|\int_0^t \sqrt{2}dW_s^{i,N}\right|  + \sup_{t\in[0,T]}\left|\int_0^t \sigma_sdB_s\right| \geq \widetilde{C}_{K_0,V}N^\beta \right\}
    \nonumber\\
    &\subset L_1^N \cup L_2^N \cup L_3^N, \label{time dec}
    \end{align}
    with
    \begin{align}
        L_1^N \doteq \left\{\exists i ;|X_0^{i,N}| \geq C_{K_0,V} N^\beta\right\},
    \end{align}
    \begin{align}
        L_2^N \doteq \left\{\exists i; \sup_{t\in [0,T]} \left|\int_0^t dW_s^{i,N}\right|  \geq C_{K_0,V} N^\beta\right\} \label{L2}
    \end{align}
    and
    \begin{align}
    L_3^N \doteq  \left\{\exists i;  \sup_{t\in[0,T]}\left|\int_0^t \sigma_sdB_s\right| \geq C_{K_0,V} N^\beta\right\} \label{predecay}
    \end{align}
with $C_{K_0, V} \doteq (1-T\|\nabla V\|_1\|K_0\|_\infty)/3\sqrt{2} >0$, for $T < (1/\|\nabla V\|_1\|K_0\|_\infty)$. 

We note that, with the same computations like those led (\ref{decay initial}), we deduce
\begin{align}    \mathbb{P}\left(L_1^N\right)&\overset{}{\leq} N(C_{K_0,V}N^\beta)^{d-2} e^{- (C_{K_0,V}N^\beta)^2} \label{decay initial1}
\end{align}

We now take care about the $L_2^N$.
First, we observe that, in view of Definition 5.1 in Chapter 2 of \cite{Karatzas}, 
\begin{align}
    \Omega = \prod_{j=1}^d \Omega^j, \, \, \, \mathcal{F} = \bigotimes_{j=1}^d \mathcal{F}^j, \, \, \, \mathbb{P}=\bigotimes_{j=1}^d \mathbb{P}^j \label{produ}
\end{align}
and then
\begin{align}
    W_t^{i,N}(\omega) = (W_t^{i,N,j}(\omega_j))_{j=1}^d
\end{align}
for $W^{i,N,j}$ a standard, one-dimensional Brownian motion on $(\Omega^j,\mathcal{F}^j,\mathbb{P}^j)$, $j \in \{1,...,d\}$.

Therefore, we note that $ W_t^{i,N} = \int_0^t  dW_s^{i,N}$, and since
\begin{align*}
    |W_t^{i,N}|^2 = \sum_{j=1}^d (W_t^{i,N,j})^2,
\end{align*}
if \(\sup_{t\leq T}|W_t^{i,N}| \geq a\), then there exists $j\in \{1,...,d\}$, such that
\(
\sup_{t\leq T}|W_t^{i,N,j}| \geq \frac{a}{\sqrt{d}}.
\)

Hence,
\begin{align}
    \left\{\sup_{t\le T} |W_t^{i,N}| \geq a\right\}
\subset
\bigcup_{j=1}^d
\left[\Omega^1 \times...\times\widetilde{\Omega}^j \times... \times \Omega^d \right]\label{d dimens}
\end{align}
with
\begin{align*}
    \widetilde{\Omega}^j \doteq \left\{
\sup_{t\le T} |W_t^{i,N,j}|
\geq
\frac{a}{\sqrt{d}}
\right\}
\end{align*}
and then, we obtain
\begin{align}
    \mathbb{P}\left(L_2^N \right) &\overset{(\ref{L2}), (\ref{d dimens})}{\leq} \sum_{i=1}^N \sum_{j=1}^d \mathbb{P}^j\left(\sup_{t\leq T} |W_t^{i,N,j}| \geq \frac{C_{K_0,V}N^\beta}{\sqrt{d}}\right) \nonumber\\
    &\leq \left(\sqrt{\frac{T}{2\pi}}\frac{4dN}{\frac{C_{K_0,V}N^\beta}{\sqrt{d}}}\right) \exp{\left(-\left(\frac{C_{K_0,V}N^\beta}{\sqrt{d}}\right)^2/(2T)\right)}, \label{decay brownian}
\end{align}
where in the second inequality, we are using (8.3)' in Chapter 2 of \cite{Karatzas}, page 96.

Now we will derive an analogous estimate for $L_3^N$. First, we observe that, if $\sigma=0$ we are done, thus suppose $\sigma \neq0$. Define,
\begin{align}
    Y_t \doteq \int_0^t\sigma_s \, dB_s,
\end{align}
which implies
\begin{align}
    Y_t^j = \sum_{k=1}^d \underbrace{\int_0^t\sigma_s^{j,k} \, dB_s^k}_{\doteq Y_s^{j,k}}.
\end{align}
Since $\sigma:[0,T]\to \mathbb{R}^d\times \mathbb{R}^d$ is bounded, for $j \in \{1,...,d\}$, the function defined by $t \mapsto \langle Y^{j,k} \rangle_t$ is continuous and increasing.

Hence, by Theorem 4.6, in Chapter 3 of \cite{Karatzas}, for all $j \in \{1,...,d\}$, we get
\begin{align}
    \sup_{t\leq T}|Y_t^j| &\leq \sum_{k=1}^d \sup_{t\leq T}\left|\int_0^t \sigma_s^{j,k} \, dB_s^k\right| \leq \sum_{k=1}^d\sup_{t\leq T}\left|\widetilde{B}^{j,k}_{\langle Y^{j,k} \rangle_t}\right| \nonumber\\
    &\leq \sum_{k=1}^d \sup_{t\leq \|\sigma\|_\infty^2 T}|\widetilde{B}^{j,k}_t|, \label{time change}
\end{align}
where $(\widetilde{B}^{j,k}_t)_{t\geq0}$ is the time-change process, which in particular, is a one dimensional standard Brownian motion.
In addition, since
\begin{align*}
    |Y_t|^2 = \sum_{j=1}^d |Y_t^j|^2
\end{align*}
if \(\sup_{t\leq T}|Y_t| \geq a\), then there exists $j \in \{1,...,d\}$, such that
\(
\sup_{t\leq T}|Y_t^j| \geq \frac{a}{\sqrt{d}}.
\)
Therefore, by (\ref{produ}) and  for all $a>0$, we find
\begin{align}
    \left\{\sup_{t\le T} |Y_t| \geq a\right\}
&\subset
\bigcup_{j=1}^d
\left\{
\sup_{t\le T} |Y_t^j|
\geq
\frac{a}{\sqrt{d}}
\right\}\nonumber\\
&\overset{(\ref{time change})}{\subset} \bigcup_{j=1}^d\left\{
\sum_{k=1}^d\sup_{t\le \|\sigma\|_\infty^2 T} |\widetilde{B}^{j,k}_t|
\geq
\frac{a}{\sqrt{d}}
\right\} \nonumber\\
&\overset{}{\subset} \bigcup_{j=1}^d \bigcup_{k=1}^d
\left[\Omega^1\times...\times \widetilde{\Omega}^{j,k} \times... \times \Omega^d\right], \label{d dimens1}
\end{align}
with
\begin{align}
    \widetilde{\Omega}^{j,k} \doteq 
\left\{
\sup_{t\le \|\sigma\|_\infty^2 T} |\widetilde{B}^{j,k}_t|
\geq
\frac{a}{d^{3/2}}
\right\}.
\end{align}
It follows that, by (\ref{time change}) and (\ref{d dimens1})
\begin{align}
    \mathbb{P}\left(L_3^N \right) &\overset{(\ref{predecay}), (\ref{d dimens1})}{\leq}  \sum_{i=1}^N\sum_{j,k=1}^d \mathbb{P}^k\left(\sup_{t\leq \|\sigma\|_\infty^2 T} |\widetilde{B}_t^{j,k}| \geq \frac{C_{K_0,V}N^\beta}{d^{3/2}}\right) \nonumber\\
    &\leq \left(\sqrt{\frac{\|\sigma\|_\infty^2 T}{2\pi}}\frac{4d^2N}{\frac{C_{K_0,V}N^\beta}{d^{3/2}}}\right) \exp{\left(-\left(\frac{C_{K_0,V}N^\beta}{d^{3/2}}\right)^2/(2\|\sigma\|_\infty^2 T)\right)}, \label{decay brownian1}
\end{align}
where in the second inequality, we are using (8.3)' in Chapter 2 of \cite{Karatzas}, page 96.

Hence, by (\ref{decay initial1}), (\ref{decay brownian}) and (\ref{decay brownian1})  in (\ref{time dec}), we find 
\begin{align}
    \sum_{N\in \mathbb{N}}\mathbb{P}\left(\tau^N < T\right) < \infty
\end{align}
and then by Borel-Cantelli Lemma, we arrive at
\begin{align}    \mathbb{P}\left(\bigcap_{M=1}^{\infty} \bigcup_{N=M}^\infty \{\tau^N < T\}\right) = 0. \label{tauN borel}
\end{align}

Finally, by  (\ref{regularityW}), (\ref{initial entropy}), (\ref{initial entropy1}),  (\ref{initial energy}), (\ref{medida total}), (\ref{marga}),  (\ref{marga1})   and (\ref{tauN borel}) in (\ref{quasienergygro}) , for $T < \min{\left(T_1;(1/\|\nabla V\|_1\|K_0\|_\infty)\right)}$
for all $\omega \in \Lambda$, with
\begin{align}
   \Lambda \doteq \Lambda_0 \cap \left(\bigcup_{M=1}^{\infty} \bigcap_{N=M}^\infty \{\tau^N = T\}\right) \cap  \left\{A_0 < \infty\right\}, \label{full meas}
\end{align}
there exists $N_0(\omega)$ such that for $N\geq N_0(\omega)$, we have $\tau^N(\omega)=T$,   
$$\mathcal{H}(\rho_0^N|\rho_0) + \|\rho_0 -\rho_0^N\|_2^2 + C_N < \exp{(-CT)}(2CT)^{-1/2}$$
and 
\begin{align}
    &\mathcal{H}(\rho_{t}^N|\rho_{t}) + \left\| \rho_{t} - \rho^N_{t}\right\|_{2}^{2} + \int_0^t \left\|\nabla(\rho_s - \rho^N_s)\right\|_{2}^{2}\, ds \nonumber\\
    &\hspace{20px}\leq  \mathcal{H}(\rho_0^N|\rho_0) + \|\rho_0 -\rho_0^N\|_2^2 + C_N \nonumber\\    &\hspace{20px}+\int_0^{t}C\left(\mathcal{H} (\rho_s^N|\rho_s) + \|\rho_s - \rho_s^N\|_2^2 + \left(\int_0^s \left\|\nabla(\rho_r - \rho^N_r)\right\|_{2}^{2}\, dr\right)\right)\, ds \nonumber\\
    &\hspace{20px} + \int_0^{t}C\left(\mathcal{H}(\rho_s^N|\rho_s) + \|\rho_s - \rho_s^N\|_2^2 + \left(\int_0^s \left\|\nabla(\rho_r - \rho^N_r)\right\|_{2}^{2}\, dr\right)\right)^3  \, ds. \label{quasienergygroh1}
     \end{align}
     In particular by the nonlinear Gronwall inequality (\ref{grownall nonlinear}) in (\ref{quasienergygroh1}), for all $t \in [0,T]$, we deduce      
\begin{align}
    &\mathcal{H}(\rho_{t}^N|\rho_{t}) + \left\| \rho_{t} - \rho^N_{t}\right\|_{2}^{2} + \int_0^t \left\|\nabla(\rho_s - \rho^N_s)\right\|_{2}^{2}\, ds \nonumber\\
    &\hspace{80px}\leq  \mathcal{H}(\rho_0^N|\rho_0) + \|\rho_0 -\rho_0^N\|_2^2 + C_N. 
     \end{align}     
     and thus, we end up with 
     
     \begin{align}
    \lim_{N \to \infty}N^{\widetilde{\theta}}\left(\sup_{t\in [0,T]}\left(\mathcal{H}(\rho_{t}^N|\rho_{t}) + \left\| \rho_{t}^N - \rho_{t}\right\|_{2}^{2}\right) + \int_0^T \left\|\nabla(\rho_s^N - \rho_s)\right\|_{2}^{2}\, ds \right)&=0, 
     \end{align}
     almost surely.
\section{Proof of Theorem \ref{first main}} \label{Proof2}
In this section we present the proof of Theorem \ref{first main}. In view of the previous section we only need the address the nonlinearity, for $d\geq 2$.
\subsection{Estimates for $I_t$: the nonlinearity} \label{nonlinea}

We now address the nonlinear terms in (\ref{ito at entropy}). The main novelty is to use Donsker-Varadhan inequality (\ref{donsker-vara}), combined with the notion of solution of limiting process given by Definition \ref{existence sol}, in the context of moderate interactions, along with the  dissipation that comes from of Fisher information and Assumption $(\mathbf{A}^{\nabla\cdot K})$,  to tackle with the singularity of kernel $K$.

First, we set
\begin{align}
    I_t &\overset{}{=}
     \int_0^t\int_{\mathbb{R}^d} \rho_s^N\nabla \ln\left(\frac{\rho_s}{\rho_s^N} \right)[K\ast \rho_s -  K\ast \rho_s^N] \,dx \,ds  \nonumber\\
     &-
     \int_0^t\int_{\mathbb{R}^d} \nabla \ln\left(\frac{\rho_s}{\rho_s^N} \right) \underbrace{\left\langle S_s^N, V^N(x - \cdot)[K\ast \rho_s^N(\cdot) - K\ast \rho_s^N(x)]\right\rangle}_{\doteq I^N} \,dx \,ds  \nonumber\\
     &\doteq I^1_t - I^2_t, \label{only Holder1}
\end{align}
where
\begin{align}
    I^1_t&= -\int_0^t\int_{\mathbb{R}^d} \left[\rho_s^N \frac{\nabla^2\rho_s}{\rho_s}  \right][ K_0\ast \rho_s -   K_0\ast \rho_s^N] \,dx \,ds \nonumber\\
    & -\int_0^t\int_{\mathbb{R}^d} \left[ \nabla \rho_s \nabla \left(\frac{\rho_s^N}{\rho_s}\right) \right][ K_0\ast \rho_s -  K_0\ast \rho_s^N] \,dx \,ds \nonumber\\
    &\doteq - I^{1,1}_t - I^{1,2}_t. \label{It1,1,2}
\end{align}
Now by $\epsilon$-Young inequality and convolution inequality, we deduce
\begin{align}
    &I_t^{1,2} =\int_0^t\int_{\mathbb{R}^d} \left[ \nabla \rho_s \left(\frac{\sqrt{\rho_s^N}}{\rho_s}\right) \left(\frac{\rho_s}{\sqrt{\rho_s^N}}\right) \nabla \left(\frac{\rho_s^N}{\rho_s}\right) \right][ K_0\ast \rho_s -  K_0\ast \rho_s^N] \,dx \,ds \nonumber \\
    &\leq \epsilon \int_0^t \int_{\mathbb{R}^d} \frac{\rho_s^2}{\rho_s^N}\left|\nabla \left(\frac{\rho_s^N}{\rho_s}\right)\right|^2 \, dx \, ds + C_{\epsilon}\int_0^t \int_{\mathbb{R}^d} \frac{|\nabla \rho_s|^2}{\rho_s^2}\rho_s^N\left|K_0 \ast(\rho_s - \rho_s^N)\right|^2 \, dx \, ds  \nonumber \\
     &\leq \epsilon \int_0^t \int_{\mathbb{R}^d} \frac{\rho_s^2}{\rho_s^N}\left|\nabla \left(\frac{\rho_s^N}{\rho_s}\right)\right|^2 dx \, ds + C_{\epsilon,K_0}\int_0^t \|\rho_s - \rho_s^N\|_1^2\int_{\mathbb{R}^d} \frac{|\nabla \rho_s|^2}{\rho_s^2}\rho_s^N dx \, ds. \label{nonlinear}
    \end{align}
   We observe that, for  $\eta > 2C_1(8T + C_3) + 1$
    \begin{align}        \int_{\mathbb{R}^d}\frac{C_3}{1+t}\exp{\left(-\frac{|x - X_t|^2}{8t + C_3}\right)}\exp{\left(\frac{C_1(1 + |x - X_t|^2)}{\eta}\right)}\, dx < \infty \label{logfinite}
    \end{align}
     which implies by Donsker-Varadhan inequality (\ref{donsker-vara}),
    \begin{align}
    \int_{\mathbb{R}^d} \rho_s^N  \frac{\left|\nabla \rho_s\right|^2}{\rho_s^2} \,dx &\overset{(\ref{donsker-vara})}{\leq} \eta \mathcal{H}(\rho_s^N|\rho_s) + \eta \ln{\left(\int_{\mathbb{R}^d}\rho_s\exp{\left(\frac{|\nabla \ln{\rho_s}|^2}{\eta}\right)}\, dx\right)} \nonumber \\
    &\overset{(\ref{logfinite}),(\ref{decay grad ln rho}),(\ref{decay rho})}{\leq} \eta \mathcal{H}(\rho_s^N|\rho_s) + C_\eta. \label{grad donsker1}
\end{align}
Hence by Csiszár-Kullback-Pinsker inequality (\ref{pinsker}) and (\ref{nonlinear}), we derive
\begin{align}
    I_t^{1,2} 
     &\leq \epsilon \int_0^t \mathcal{I}(\rho^N_s|\rho_s) \, ds + C_{\epsilon,\eta,K_0}\int_0^t \mathcal{H}(\rho^N_s|\rho_s) \,ds. \label{nonlinear1}
    \end{align}
Regarding term $I_t^{1,1}$ in (\ref{It1,1,2}), by convolution inequality, we find
\begin{align} 
    I_t^{1,1}
    &\leq\int_0^t\|K_0\|_\infty\|\rho^N_s - \rho_s\|_1\int_{\mathbb{R}^d}  \left|\frac{\nabla^2\rho_s}{\rho_s}  \right| \rho^N_s \,dx \,ds. \label{nonlinear2}
    \end{align}
    In addition, if $C_{1,2} \doteq \max{(C_1,C_2)}$ since
    \begin{align}
        \left|\frac{\nabla^2 \rho_s(x)}{\rho_s(x)}\right| \leq |\nabla^2 \ln{(\rho_s)}(x)| + |\nabla \ln{(\rho_s}(x))|^2 \overset{(\ref{decay grad ln rho}),(\ref{decay hess ln rho})}{\leq} C_{1,2}(1 + |x - X_t|^2), \label{sum decay}
    \end{align}
   for  $\eta > 2C_{1,2}(8T + C_3) + 1$, we obtain
    \begin{align}        \int_{\mathbb{R}^d}\frac{C_3}{1+t}\exp{\left(-\frac{|x - X_t|^2}{8t + C_3}\right)}\exp{\left(\frac{C_{1,2}(1 + |x - X_t|^2)}{\eta}\right)}\, dx < \infty. \label{logfinite1}
    \end{align}
    It follows by  Donsker-Varadhan inequality (\ref{donsker-vara}), 
    \begin{align}
    \int_{\mathbb{R}^d} \rho_s^N  \left|\frac{\nabla^2 \rho_s}{\rho_s} \right|\,dx &\overset{(\ref{donsker-vara})}{\leq} \eta \mathcal{H}(\rho_s^N|\rho_s) + \eta \ln{\left(\int_{\mathbb{R}^d}\rho_s\exp{\left(\frac{\left|\frac{\nabla^2 \rho_s}{ \rho_s}\right|}{\eta}\right)}\, dx\right)}  \nonumber\\
    &\overset{(\ref{logfinite1}),(\ref{decay rho})}{\leq}  \eta \mathcal{H}(\rho_s^N|\rho_s) + C_\eta. \label{donsker}
\end{align}

Therefore, by (\ref{donsker}) in (\ref{nonlinear2}),  we get 
    \begin{align}
     I_t^{1,1} &\overset{}{\leq} C_{\eta,K_0}\int_0^t \mathcal{H}(\rho_s^N|\rho_s)\, ds + TC_\eta.
   \label{It11}
\end{align}
Thus, from  $(\ref{nonlinear1})$ and $(\ref{It11})$, we have  the following estimate for the term $I_t^1$ in (\ref{only Holder}):
\begin{align}
   I_t^1&\leq C_{\epsilon,K_0,\eta}\int_0^t \mathcal{H}(\rho_s^N|\rho_s)\, ds + \epsilon  \int_0^t \mathcal{I}(\rho_s^N|\rho_s) \,ds + TC_\eta.  \label{It1}
\end{align}
Now, we will focus on the term $I_t^2$ in (\ref{only Holder}). 
First, we recall that by Assumption $\left(\mathbf{A}^{\nabla \cdot K}\right)$, $K = \nabla\cdot K_0$ with $K_0 \in L^{\infty}$.
Hence, we arrive at
\begin{align}
    \left|K \ast \rho_s^N(\cdot) - K\ast \rho_s^N(x)\right| & = \left|K_0 \ast \nabla \rho_s^N(\cdot) - K_0\ast \nabla \rho_s^N(x)\right|\nonumber\\
    &\leq \int_{\mathbb{R}^d} \left|K_0(\cdot -y) - K_0(x - y)\right||\nabla \rho_s^N(y)| \, dy\nonumber\\
    &\leq \int_{\mathbb{R}^d} 2\|K_0\|_{\infty}|\nabla\rho_s^N(y)| \, dy \nonumber\\
    &= 2\|K_0\|_{\infty} \int_{\mathbb{R}^d} |\nabla\rho_s^N(y)| \, dy. \label{K_0 Holder1}
\end{align}
Additionally,  $\rho_s^N \in \mathcal{P}\left(\mathbb{R}^d\right)$, $\mathbb{P}$-a.s. and then by Holder's inequality we obtain 
\begin{align}
   \int_{\mathbb{R}^d}|\nabla \rho_s^N| \,dx &= \int_{\mathbb{R}^d}|\nabla \rho_s^N| \frac{\sqrt{\rho_s^N}}{\sqrt{\rho_s^N}} \,dx \nonumber\\
    &\overset{\text{Holder}}{\leq} \left(\int_{\mathbb{R}^d}\frac{|\nabla \rho_s^N|^2}{\rho_s^N} \,dx\right)^{\frac{1}{2}} \left(\int_{\mathbb{R}^d}\rho_s^N \,dx \right)^{\frac{1}{2}} \nonumber \\
    &\overset{}{=} \left(\int_{\mathbb{R}^d}\frac{|\nabla \rho_s^N|^2}{\rho_s^N}\,dx \right)^{\frac{1}{2}}. \label{fisher and L11}
\end{align}
Now from  Leibniz's rule we deduce, 
\begin{align} \label{fisher rho^N and difference1}
   \int_{\mathbb{R}^d}  \frac{|\nabla \rho_s^N|^2}{\rho_s^N} \,dx \nonumber& = \int_{\mathbb{R}^d}  \frac{\left|\nabla \left(\rho_s\frac{\rho_s^N}{\rho_s}\right)\right|^2}{\rho_s^N} \,dx\\
    \nonumber&\leq2\int_{\mathbb{R}^d}  \frac{\left|\nabla \rho_s\left(\frac{\rho_s^N}{\rho_s}\right)\right|^2}{\rho_s^N} \,dx + 2\int_{\mathbb{R}^d}  \frac{\left|\rho_s\nabla \left(\frac{\rho_s^N}{\rho_s}\right)\right|^2}{\rho_s^N} \,dx
    \\
   &\overset{}{\leq} 2\int_{\mathbb{R}^d} \rho_s^N  \frac{\left|\nabla \rho_s\right|^2}{\rho_s^2} \,dx + 2\int_{\mathbb{R}^d} \rho_s^N \left|\nabla \ln \left(\frac{\rho_s^N}{\rho_s}\right)\right|^2 \,dx
    \nonumber\\
    &\overset{(\ref{logfinite1})}{\leq} 2\eta \mathcal{H}(\rho_s^N|\rho_s)  + C_\eta + 2\int_{\mathbb{R}^d} \rho_s^N \left|\nabla \ln \left(\frac{\rho_s^N}{\rho_s}\right)\right|^2 \,dx. 
   \end{align}

From $(\ref{K_0 Holder1})$, $(\ref{fisher and L11})$ and $(\ref{fisher rho^N and difference1})$, recalling $I_t^2$ in (\ref{only Holder}), since by Assumption $\left(\mathbf{A}^{\nabla\cdot K}\right)$ $\|K_0\|_{\infty}\leq 1/4$, we have 
\begin{align*}
   I^N
   &\leq 2\|K_0\|_{\infty}\left\langle S_s^N, V^N(x - \cdot)\right\rangle \left(2\eta \mathcal{H}(\rho^N_s|\rho_s) + 2\mathcal{I}\left(\rho_s^N|\rho_s\right) + C_\eta \right)^{\frac{1}{2}} \nonumber\\
&\leq 1/2\rho^N_s  \left(2\eta \mathcal{H}(\rho^N_s|\rho_s) + 2\mathcal{I}\left(\rho_s^N|\rho_s\right) + C_\eta \right)^{\frac{1}{2}},
\end{align*}
and then
\begin{align}
    I^{2}_t &\doteq \int_0^t\int_{\mathbb{R}^d} \nabla \ln\left(\frac{\rho_s}{\rho_s^N} \right) \left\langle S_s^N, V^N(x - \cdot)[K\ast \rho_s^N(\cdot) - K\ast \rho_s^N(x)]\right\rangle \,dx \,ds \nonumber\\
    &\leq \int_0^t\int_{\mathbb{R}^d} \left| \nabla \ln\left(\frac{\rho_s}{\rho_s^N} \right)\right|(1/2) \rho^N_s \left(2\eta \mathcal{H}(\rho^N_s|\rho_s) + 2\mathcal{I}\left(\rho_s^N|\rho_s\right) + C_\eta\right)^{\frac{1}{2}}\,dx \,ds \nonumber \\
    &\leq \int_0^t\int_{\mathbb{R}^d} \left| \nabla \ln\left(\frac{\rho_s}{\rho_s^N} \right)\right|(1/2) (\rho^N_s)^{1/2} (\rho^N_s)^{1/2}\left(2\eta \mathcal{H}(\rho^N_s|\rho_s)\right)^{\frac{1}{2}} \,dx \,ds  \nonumber\\
    &+ \int_0^t\int_{\mathbb{R}^d} \left| \nabla \ln\left(\frac{\rho_s}{\rho_s^N} \right)\right|(1/2) \rho^N_s  \left(2\mathcal{I}\left(\rho_s^N|\rho_s\right)\right)^{\frac{1}{2}}\,dx \,ds \nonumber \\
    &+ \int_0^t\int_{\mathbb{R}^d} \left| \nabla \ln\left(\frac{\rho_s}{\rho_s^N} \right)\right|(1/2) (\rho^N_s)^{1/2}(\rho^N_s)^{1/2}  \left(C_\eta\right)^{\frac{1}{2}}\,dx \,ds. 
    \end{align}
    It follows that, by Young and Holder inequalities,
    \begin{align}
    &I_t^2\overset{\text{Young}}{\leq} \epsilon \int_0^t\int_{\mathbb{R}^d} \rho^N_s \left| \nabla \ln\left(\frac{\rho_s}{\rho_s^N} \right)\right|^2 \,dx \,ds + C_{\epsilon,\eta}\int_0^t\ \mathcal{H}(\rho^N_s|\rho_s) \,ds  \nonumber\\
    &+ \sqrt{2}/2\int_0^t\left(\mathcal{I}\left(\rho_s^N|\rho_s\right)\right)^{\frac{1}{2}}\int_{\mathbb{R}^d} \left| \nabla \ln\left(\frac{\rho_s}{\rho_s^N} \right)\right| (\rho^N_s)^{1/2}(\rho^N_s)^{1/2}  \,dx \,ds \nonumber \\
    &+ \epsilon \int_0^t\int_{\mathbb{R}^d} \rho^N_s\left| \nabla \ln\left(\frac{\rho_s}{\rho_s^N} \right)\right|^2 \,dx \,ds + TC_{\epsilon,\eta} \nonumber \\
    &\overset{\text{Holder}}{\leq} (2\epsilon +\sqrt{2}/2)\int_0^t\mathcal{I}\left(\rho_s^N|\rho_s\right)\,ds +  C_{\epsilon,\eta}\int_0^t\ \mathcal{H}(\rho^N_s|\rho_s) \,ds + TC_{\epsilon,\eta}. \label{no Holder1}    \end{align}

Finally, by $(\ref{It1})$ and $(\ref{no Holder1})$ in $(\ref{only Holder1})$, for  $\epsilon <<1$ and $\eta >>1$,  we conclude  with the following estimate for the nonlinear term:
\begin{align}
   I_t
   &\leq C_{\epsilon,K_0,\eta} \int_0^t \mathcal{H}(\rho_s^N|\rho_s)  \,ds+ (3\epsilon + \sqrt{2}/2) \int_0^t \mathcal{I}(\rho_s^N|\rho_s) \,ds + TC_{\epsilon,\eta}.  \label{I_t regular kernel}
\end{align}
 \subsection{Application of Gronwall's Lemma
}\label{gronwall}
In this subsection, we conclude our localizated estimates, by applying
Gronwall's inequality.
Indeed, we put $(\ref{II_t}), (\ref{cancelations}),(\ref{I_t regular kernel})$ and (\ref{marga}) into $(\ref{ito at entropy})$, to get
\begin{align}
    \mathcal{H}(\rho_{t\wedge \tau^N}^N|\rho_{t\wedge \tau^N}) - \mathcal{H}(\rho_0^N|\rho_0) &\leq C_{\epsilon,K_0,\eta}\int_0^{t\wedge \tau^N}\mathcal{H}(\rho_s^N|\rho_s)\, ds \nonumber\\
    &\overset{}{+} (3\epsilon + \sqrt{2}/2)\int_0^{t\wedge \tau^N}\mathcal{I}(\rho_s^N|\rho_s) \,ds\nonumber\\
    &-\int_0^{t\wedge \tau^N} \mathcal{I}(\rho_s^N|\rho_s)  \,ds \nonumber\\
    &+C_{\alpha,T} N^{-\theta_1} + A_0N^{-\theta_2 + \delta}  + TC_\eta
    \label{entropy regular grownall1}
\end{align}
 and then for $\eta >>1$ and  $\epsilon << 1$,  we deduce
\begin{align}
    \mathcal{H}(\rho_{t\wedge \tau^N}^N|\rho_{t\wedge \tau^N}) - \mathcal{H}(\rho_0^N|\rho_0) &\lesssim \int_0^{t\wedge \tau^N}\mathcal{H}(\rho_s^N|\rho_s)\, ds + N^{-\theta_1}+A_0N^{-\theta_2 + \delta} + T, \label{entropy regular grownall without expectation 1}
\end{align}
with  $\theta_1 \doteq 1 - \beta(1 +\frac{2}{d} + 2\alpha)$ and $\theta_2 \doteq \left(\frac{1}{2}-\beta\Big(1 + \frac{1 }{d}\Big) \right)$.
\vspace{.2cm}

Finally, from  Gronwall's inequality (\ref{lineargronwall}) in (\ref{entropy regular grownall without expectation 1}), with the analogous computations as led  (\ref{full meas}), there exists $\Lambda \in \mathcal{F}$, $\mathbb{P}(\Lambda) = 1$, with the property that,
for all $\omega \in \Lambda$,  there exists $N_0(\omega)$ such that for $N\geq N_0(\omega)$, $\tau^N(\omega)= T$ and

\begin{align}
   \sup_{t \in [0,T]}\mathcal{H}(\rho_{t\wedge \tau^N}^N|\rho_{t}) &\lesssim  \left(\mathcal{H}(\rho_0^N|\rho_0) + N^{-\theta} + A_0N^{-\theta} + T \right) \label{grownall}
   \end{align}
   with
     \begin{align*}
\theta \doteq \min \left(1 - \beta(1 +\frac{2}{d} + 2\alpha);\left(\frac{1}{2}-\beta\Big(1 + \frac{1 }{d}\Big) \right) - \delta  \right).
\end{align*}
 In particular, we get
   \begin{align}
  \limsup_{N\to \infty} \sup_{t \in [0,T]}\mathcal{H}(\rho_{t}^N|\rho_{t}) &\lesssim  T, \, \, \, \, \, \, \mathbb{P}-a.s. \label{grownallffffgh}
   \end{align}
\section{Appendix} \label{ap}
This appendix contains the proof of Theorem \ref{SPDE_Ito} together with several auxiliary inequalities used in the previous sections. We begin by establishing Theorem \ref{SPDE_Ito}, regarding the existence of a suitable solution of the limiting equation.

\subsection{Proof of Theorem \ref{existence}}

When $\sigma=0$, the equation  (\ref{SPDE_Ito}) reads
\begin{align}
     \partial_t \Tilde{\rho} - \Delta \Tilde{\rho} + \nabla \cdot (\Tilde{\rho}_tK\ast \Tilde{\rho}_t)=0, \label{sigma0}
\end{align}
and a solution to (\ref{sigma0}), in the sense of Theorem \ref{existence},  was constructed in sections 3 and 4 of \cite{Feng} and  \cite{Feng1}.
Now, taking into account (\ref{random}), we define the following process: $\rho(t,x) \doteq \Tilde{\rho}(t,x - X_t)$.

Hence, by applying the Itô’s formula, since that $\sigma$ does not depends on spatial variable, we find that $\rho$ solves (\ref{SPDE_Ito}), in the sense of Theorem \ref{existence}, for $\sigma\neq0$.

\subsection{Inequalities}
In this section, we present several inequalities used throughout the text. The first one is the classical Csiszár-Kullback-Pinsker inequality, (see \cite{Gui}) which is a fundamental result in information theory. It quantifies how control of the relative entropy implies control of the total variation distance.
\begin{lem}[Csiszár-Kullback-Pinsker inequality] \label{L1entropy}
    It holds for $f,g \in \mathcal{P}\left(\mathbb{R}^d\right)$, 
    \begin{align}
     \|f - g\|_1^2\lesssim  \mathcal{H}(f|g). \label{pinsker}
    \end{align}
\end{lem}

Now, we recall the Donsker–Varadhan inequality (see \cite{Jabin_2018}) which follows from the variational representation of the relative entropy. It provides a tool to estimate the expectation of a test function under a non-factorized probability density and, in particular, allows for a change of measure between two probability measures.
\begin{lem}[Donsker-Varadhan Inequality]
    Let $f,g \in \mathcal{P}\left(\mathbb{R}^d\right)$, $\eta >0$ and a function $\Phi$ on $\mathbb{R}^d$ such that $g\exp{\left(\Phi/\eta\right)} \in L^1$. Then, we have
    \begin{align}
        \int_{\mathbb{R}^d} f\Phi\, dx\leq \eta \mathcal{H}\left(f|g\right) + \eta\ln{\left(\int_{\mathbb{R}^d} g\exp{\left(\Phi/\eta\right)}\, dx\right)}. \label{donsker-vara}
    \end{align}
\end{lem}
The next result is a nonlinear Gronwall's inequality. A proof is provided in Theorem 21 of \cite{gronwal}.
\begin{lem} \label{nolinearg}
Let \(u\) be a nonnegative function satisfying
\[
u(t) \le c(t) + \int_{0}^{t} \big( a(s)\,u(s) + b(s)\,u^{\alpha}(s) \big)\,ds,
\qquad c(t) > 0,\ \alpha \ge 0,
\]
where \(a(t)\), \(b(t)\) and \(c(t)\) are integrable and nonnegative functions on $[0,T]$.

\medskip
\noindent
\noindent
For \(\alpha = 1\),
\begin{align}
    u(t) \le \left(\sup_{t \in [0,T]} c(t) \right) \exp\!\left( \int_{0}^{t} \big(a(s)+b(s)\big)\,ds \right). \label{lineargronwall}
\end{align}
\noindent
For \(\alpha > 1\) and $c$ a constant function, with the following additional hypothesis, for some $h >0$,
\begin{align}
    c < \left\{
\exp\!\left[(1-\alpha)\int_{0}^{h} a(s)\,ds\right]
\right\}^{\frac{1}{\alpha - 1}}
\left\{(\alpha-1)\int_{0}^{h} b(s)\,ds\right\}^{\frac{1}{1-\alpha}},
\end{align}
we get, for $t \in [0,h]$,
\begin{equation}
\begin{aligned}
u(t) \le c \Biggl\{ &\exp \left[ (1-\alpha) \int_{0}^{t} a(s) \, ds \right] \\
&- c^{-1} (\alpha - 1) \int_{0}^{t} b(s) \exp \left[ (1-\alpha) \int_{s}^{t} a(r) \, dr \right] ds \Biggr\}^{\frac{1}{\alpha-1}}. \label{grownall nonlinear}
\end{aligned}
\end{equation}
\end{lem}

\section*{Acknowledgements}
 A. B. de Souza is partially supported by  Coordenação de Aperfeiçoamento de Pessoal de Nível Superior – Brasil (CAPES) – Finance Code $001$.

\end{document}